\documentclass[secthm,seceqn,amsthm,ussrhead,12pt]{amsart}
\usepackage[utf8]{inputenc}
\usepackage[english]{babel}
\usepackage{amssymb,amsmath,amsthm,amsfonts,xcolor,enumerate,hyperref,comment,longtable,cleveref}

\usepackage{times}
\usepackage{cite}
\usepackage{pdflscape}
\usepackage{ulem}
\usepackage[mathcal]{euscript}
\usepackage{tikz}
\usepackage{hyperref}
\usepackage{cancel}
\usepackage{stmaryrd}

\usetikzlibrary{arrows}

\newtheorem*{theoremA}{Theorem A}
\newtheorem*{theoremB}{Theorem B}

\newtheorem{theorem}{Theorem}
\newtheorem{lemma}[theorem]{Lemma}

\newtheorem{remark}[theorem]{Remark}

\newtheorem{definition}[theorem]{Definition}

\usepackage{stmaryrd}
\usepackage{xcolor}

\begin{document}
\noindent{\Large 
One-generated nilpotent assosymmetric algebras}
 
\medskip

   \

   {\bf
   Ivan   Kaygorodov \& Farukh Mashurov
 }

 \

{\tiny

\smallskip
\smallskip
   E-mail addresses:

\smallskip

Ivan   Kaygorodov (kaygorodov.ivan@gmail.com)

Farukh Mashurov (f.mashurov@gmail.com)

}

\smallskip
\smallskip
\noindent{\bf Abstract}:
{\it We give the classification of $5$- and $6$-dimensional complex one-generated nilpotent assosymmetric algebras.}

\ 

\smallskip
\noindent {\bf Keywords}:
{\it assosymmetric algebras, nilpotent algebras, algebraic classification, central extension.}

\ 

\smallskip
\noindent {\bf MSC2010}: 17A30, 17D25.

\section*{Introduction}

Algebraic classification (up to isomorphism) of algebras of small dimension from a certain variety defined by a family of polynomial identities is a classic problem in the theory of non-associative algebras. There are many results related to algebraic classification of small dimensional algebras in varieties of Jordan, Lie, Leibniz, Zinbiel and other algebras.
Another interesting approach of studying algebras of a fixed dimension is to study them from a geometric point of view (that is, to study degenerations and deformations of these algebras). The results in which the complete information about degenerations of a certain variety is obtained are generally referred to as the geometric classification of the algebras of these variety. There are many results related to geometric classification of Jordan, Lie, Leibniz, Zinbiel and other algebras \cite{ack, cfk19,  ikv17, kv16}.
Another interesting direction is a study of one-generated objects.
Well know the description of one-generated finite groups: there is only one one-generated group of order $n$.
    In the case of algebras, there are some similar results,
    such that the description of $n$-dimensional one-generated nilpotent associative \cite{karel}, noncommutative Jordan \cite{jkk19},  Leibniz and Zinbiel algebras\cite{bakhrom}. 
    It was proven, that there is only one $n$-dimensional one-generated nilpotent algebra in these varieties.
    But on the other side, as we can see in varieties of Novikov \cite{kkk18}, assosymmetric \cite{ikm19},  bicommutative \cite{kpv19}, commutative \cite{fkkv}, and terminal \cite{kkp19} algebras, there are more than one $4$-dimensional one-generated nilpotent algebra from these varieties.
    One-generated nilpotent Novikov algebras in dimensions 5 and 6 were studied in \cite{ckkk19}, one-generated nilpotent terminal algebras in dimension 5 were studied in \cite{kks}.
In the present paper, we give the algebraic   classification of
$5$- and $6$-dimensional complex one-generated   nilpotent assosymmetric   algebras,
which first appeared in the paper by Kleinfeld in 1957 \cite{Kleinfeld}.

The variety of assosymmetric algebras is defined by the following identities of right- and left-symmetric:
\[
\begin{array}{rclllrcl}
(x,y,z) &=& (x,z,y), & \ &  (x,y,z) &=& (y,x,z),
\end{array} \]
where $(x,y,z)=(xy)z-x(yz).$
It admits the commutative associative and associative algebras as a subvariety.  Kleinfeld proved that an assosymmetric ring of characteristic different from 2 and 3 without ideals $I\neq0,$ such that $I^2=0$ is associative \cite{Kleinfeld}. The free base elements of assosymmetric algebras were described in \cite{hentzel}.
The algebraic and geometric classification of $4$-dimensional complex nilpotent assosymmetric algebras was given in \cite{ikm19}.
Also, assosymmetric algebras were studied in \cite{boers,boers95,dzhumadildayev, pr77, kim, askarbekzat}. 

The key step in our method for algebraically classifying assosymmetric nilpotent algebras is the calculation of central extensions of smaller algebras.  It comes as no surprise that the central extensions of Lie and non-Lie algebras have been exhaustively studied for years. It is interesting both to describe them and to use them to classify different varieties of algebras \cite{gkks,ckkk20,krs19,lisa,gkks,hac16,kkl18,ss78}.
Firstly, Skjelbred and Sund devised a method for classifying nilpotent Lie algebras employing central extensions  \cite{ss78}.
Using this method, all the non-Lie central extensions of  all $4$-dimensional Malcev algebras were described  afterwards \cite{hac16}, and also
all the anticommutative central extensions of the $3$-dimensional anticommutative algebras \cite{cfk182},
and all the central extensions of the $2$-dimensional algebras \cite{cfk18}.
Moreover, the method is especially indicated for the classification of nilpotent algebras and it was used to describe
all the $4$-dimensional nilpotent associative algebras \cite{degr1},
all the $4$-dimensional nilpotent Novikov algebras \cite{kkk18},
all the $4$-dimensional nilpotent bicommutative algebras \cite{kpv19},
all the $5$-dimensional nilpotent Jordan algebras \cite{ha16j},
all the $5$-dimensional nilpotent restricted Lie algebras \cite{usefi1}, 
all the $6$-dimensional nilpotent Lie algebras \cite{degr3,degr2},
all the $6$-dimensional nilpotent Malcev algebras \cite{hac18},
all the $6$-dimensional nilpotent Tortkara  algebras \cite{gkk18}
and some others.

\section{The algebraic classification of nilpotent assosymmetric algebras}
\subsection{Method of classification of nilpotent algebras}

The objective of this section is to give an analogue of the Skjelbred-Sund method for classifying nilpotent assosymmetric algebras. As other analogues of this method were carefully explained in, for example, \cite{hac16,cfk18}, we will give only some important definitions, and refer the interested reader to the previous sources. We will also employ their notations.

Let $({\bf A}, \cdot)$ be an assosymmetric  algebra over  $\mathbb C$
and ${\bf V}$ a vector space over ${\mathbb C}$. We define the $\mathbb C$-linear space ${\rm Z^{2}}\left(
\bf A,{\bf V} \right) $  as the set of all  bilinear maps $\theta  \colon {\bf A} \times {\bf A} \longrightarrow {{\bf V}}$
such that
\[ \theta(xy,z)-\theta(x,yz)=\theta(xz,y)-\theta(x,zy), \]
\[\theta(xy,z)-\theta(x,yz)=\theta(yx,z)-\theta(y,xz). \]
These maps will be called {\it cocycles}. Consider a
linear map $f$ from $\bf A$ to  ${\bf V}$, and set $\delta f\colon {\bf A} \times
{\bf A} \longrightarrow {{\bf V}}$ with $\delta f  (x,y ) =f(xy )$. Then, $\delta f$ is a cocycle, and we define ${\rm B^{2}}\left(
{\bf A},{{\bf V}}\right) =\left\{ \theta =\delta f\ : f\in {\rm Hom}\left( {\bf A},{{\bf V}}\right) \right\} $, a linear subspace of ${\rm Z^{2}}\left( {\bf A},{{\bf V}}\right) $; its elements are called
{\it coboundaries}. The {\it second cohomology space} ${\rm H^{2}}\left( {\bf A},{{\bf V}}\right) $ is defined to be the quotient space ${\rm Z^{2}}
\left( {\bf A},{{\bf V}}\right) \big/{\rm B^{2}}\left( {\bf A},{{\bf V}}\right) $.
 
\

Let ${\rm Aut}({\bf A}) $ be the automorphism group of the assosymmetric algebra ${\bf A} $ and let $\phi \in {\rm Aut}({\bf A})$. Every $\theta \in
{\rm {\rm Z^{2}}}\left( {\bf A},{{\bf V}}\right) $ defines $\phi \theta (x,y)
=\theta \left( \phi \left( x\right) ,\phi \left( y\right) \right) $, with $\phi \theta \in {\rm {\rm Z^{2}}}\left( {\bf A},{{\bf V}}\right) $. It is easily checked that ${\rm Aut}({\bf A})$
acts on ${\rm {\rm Z^{2}}}\left( {\bf A},{{\bf V}}\right) $, and that
 ${\rm B^{2}}\left( {\bf A},{{\bf V}}\right) $ is invariant under the action of ${\rm Aut}({\bf A}).$  
 So, we have that ${\rm Aut}({\bf A})$ acts on ${\rm H^{2}}\left( {\bf A},{{\bf V}}\right)$.

\

Let $\bf A$ be an assosymmetric  algebra of dimension $m<n$ over  $\mathbb C$, ${{\bf V}}$ a $\mathbb C$-vector
space of dimension $n-m$ and $\theta$ a cocycle, and consider the direct sum ${\bf A}_{\theta } = {\bf A}\oplus {{\bf V}}$ with the
bilinear product `` $\left[ -,-\right] _{{\bf A}_{\theta }}$'' defined by $\left[ x+x^{\prime },y+y^{\prime }\right] _{{\bf A}_{\theta }}=
 xy +\theta(x,y) $ for all $x,y\in {\bf A},x^{\prime },y^{\prime }\in {{\bf V}}$.
It is straightforward that ${\bf A_{\theta}}$ is an assosymmetric algebra if and only if $\theta \in {\rm Z}^2({\bf A}, {{\bf V}})$; it is  called an $(n-m)$-{\it dimensional central extension} of ${\bf A}$ by ${{\bf V}}$.

We also call the
set ${\rm Ann}(\theta)=\left\{ x\in {\bf A}:\theta \left( x, {\bf A} \right)+ \theta \left({\bf A} ,x\right) =0\right\} $
the {\it annihilator} of $\theta $. We recall that the {\it annihilator} of an  algebra ${\bf A}$ is defined as
the ideal ${\rm Ann}(  {\bf A} ) =\left\{ x\in {\bf A}:  x{\bf A}+ {\bf A}x =0\right\}$. Observe
 that
${\rm Ann}\left( {\bf A}_{\theta }\right) =\big({\rm Ann}(\theta) \cap{\rm Ann}({\bf A})\big)
 \oplus {{\bf V}}$.

\

\begin{definition}
Let ${\bf A}$ be an algebra and $I$ be a subspace of ${\rm Ann}({\bf A})$. If ${\bf A}={\bf A}_0 \oplus I$
then $I$ is called an {\it annihilator component} of ${\bf A}$.
A central extension of an algebra $\bf A$ without annihilator component is called a {\it non-split central extension}.
\end{definition}

\

The following result is fundamental for the classification method.

\begin{lemma}
Let ${\bf A}$ be an $n$-dimensional assosymmetric algebra such that $\dim \ {\rm Ann}({\bf A})=m\neq0$. Then there exists, up to isomorphism, a unique $(n-m)$-dimensional assosymmetric  algebra ${\bf A}'$ and a bilinear map $\theta \in {\rm Z}^2({\bf A}, {{\bf V}})$ with ${\rm Ann}({\bf A})\cap{\rm Ann}(\theta)=0$, where ${\bf V}$ is a vector space of dimension m, such that ${\bf A} \cong {{\bf A}'}_{\theta}$ and
 ${\bf A}/{\rm Ann}({\bf A})\cong {\bf A}'$.
\end{lemma}

For the proof, we refer the reader to~\cite[Lemma 5]{hac16}.

\

Now, we seek a condition on the cocycles to know when two $(n-m)$-central extensions are isomorphic.
Let us fix a basis $e_{1},\ldots ,e_{s}$ of ${{\bf V}}$, and $
\theta \in {\rm Z^{2}}\left( {\bf A},{{\bf V}}\right) $. Then $\theta $ can be uniquely
written as $\theta \left( x,y\right) =
\displaystyle \sum_{i=1}^{s} \theta _{i}\left( x,y\right) e_{i}$, where $\theta _{i}\in
{\rm Z^{2}}\left( {\bf A},\mathbb C\right) $. It holds that $\theta \in
{\rm B^{2}}\left( {\bf A},{{\bf V}}\right) $\ if and only if all $\theta _{i}\in {\rm B^{2}}\left( {\bf A},
\mathbb C\right) $, and it also holds that ${\rm Ann}(\theta)={\rm Ann}(\theta _{1})\cap{\rm Ann}(\theta _{2})\ldots \cap{\rm Ann}(\theta _{s})$. 
Furthermore, if ${\rm Ann}(\theta)\cap {\rm Ann}\left( {\bf A}\right) =0$, then ${\bf A}_{\theta }$ has an
annihilator component if and only if $\left[ \theta _{1}\right] ,\left[
\theta _{2}\right] ,\ldots ,\left[ \theta _{s}\right] $ are linearly
dependent in ${\rm H^{2}}\left( {\bf A},\mathbb C\right)$ (see \cite[Lemma 13]{hac16}).

\;

Recall that, given a finite-dimensional vector space ${{\bf V}}$ over $\mathbb C$, the {\it Grassmannian} ${\rm G}_{k}\left( {{\bf V}}\right) $ is the set of all $k$-dimensional
linear subspaces of $ {{\bf V}}$. Let ${\rm G}_{s}\left( {\rm H^{2}}\left( {\bf A},\mathbb C\right) \right) $ be the Grassmannian of subspaces of dimension $s$ in
${\rm H^{2}}\left( {\bf A},\mathbb C\right) $.
 For ${\rm W}=\left\langle
\left[ \theta _{1}\right] ,\left[ \theta _{2}\right] ,\dots,\left[ \theta _{s}
\right] \right\rangle \in {\rm G}_{s}\left( {\rm H^{2}}\left( {\bf A},\mathbb C
\right) \right) $ and $\phi \in {\rm Aut}({\bf A})$, define $\phi {\rm W}=\left\langle \left[ \phi \theta _{1}\right]
,\left[ \phi \theta _{2}\right] ,\dots,\left[ \phi \theta _{s}\right]
\right\rangle $. It holds that $\phi {\rm W}\in {\rm G}_{s}\left( {\rm H^{2}}\left( {\bf A},\mathbb C \right) \right) $, and this induces an action of ${\rm Aut}({\bf A})$ on ${\rm G}_{s}\left( {\rm H^{2}}\left( {\bf A},\mathbb C\right) \right) $. We denote the orbit of ${\rm W}\in {\rm G}_{s}\left(
{\rm H^{2}}\left( {\bf A},\mathbb C\right) \right) $ under this action  by ${\rm Orb}({\rm W})$. Let
\[
{\rm W}_{1}=\left\langle \left[ \theta _{1}\right] ,\left[ \theta _{2}\right] ,\dots,
\left[ \theta _{s}\right] \right\rangle ,{\rm W}_{2}=\left\langle \left[ \vartheta
_{1}\right] ,\left[ \vartheta _{2}\right] ,\dots,\left[ \vartheta _{s}\right]
\right\rangle \in {\rm G}_{s}\left( {\rm H^{2}}\left( {\bf A},\mathbb C\right)
\right).
\]
Similarly to~\cite[Lemma 15]{hac16}, in case ${\rm W}_{1}={\rm W}_{2}$, it holds that \[ \bigcap\limits_{i=1}^{s}{\rm Ann}(\theta _{i})\cap {\rm Ann}\left( {\bf A}\right) = \bigcap\limits_{i=1}^{s}
{\rm Ann}(\vartheta _{i})\cap{\rm Ann}( {\bf A}) ,\] 
and therefore the set
\[
{\rm T}_{s}({\bf A}) =\left\{ {\rm W}=\left\langle \left[ \theta _{1}\right] ,
\left[ \theta _{2}\right] ,\dots,\left[ \theta _{s}\right] \right\rangle \in
{\rm G}_{s}\left( {\rm H^{2}}\left( {\bf A},\mathbb C\right) \right) : \bigcap\limits_{i=1}^{s}{\rm Ann}(\theta _{i})\cap{\rm Ann}({\bf A}) =0\right\}
\]
is well defined, and it is also stable under the action of ${\rm Aut}({\bf A})$ (see~\cite[Lemma 16]{hac16}).

\

Now, let ${{\bf V}}$ be an $s$-dimensional linear space and let us denote by
${\rm E}\left( {\bf A},{{\bf V}}\right) $ the set of all non-split $s$-dimensional central extensions of ${\bf A}$ by
${{\bf V}}$. We can write
\[
{\rm E}\left( {\bf A},{{\bf V}}\right) =\left\{ {\bf A}_{\theta }:\theta \left( x,y\right) = \sum_{i=1}^{s}\theta _{i}\left( x,y\right) e_{i} \ \ \text{and} \ \ \left\langle \left[ \theta _{1}\right] ,\left[ \theta _{2}\right] ,\dots,
\left[ \theta _{s}\right] \right\rangle \in {\rm T}_{s}({\bf A}) \right\} .
\]

Finally, we are prepared to state our main result, which can be proved as \cite[Lemma 17]{hac16}.

\begin{lemma}
 Let ${\bf A}_{\theta },{\bf A}_{\vartheta }\in {\rm E}\left( {\bf A},{{\bf V}}\right) $. Suppose that $\theta \left( x,y\right) =  \displaystyle \sum_{i=1}^{s}
\theta _{i}\left( x,y\right) e_{i}$ and $\vartheta \left( x,y\right) =
\displaystyle \sum_{i=1}^{s} \vartheta _{i}\left( x,y\right) e_{i}$.
Then the assosymmetric algebras ${\bf A}_{\theta }$ and ${\bf A}_{\vartheta } $ are isomorphic
if and only if
$${\rm Orb}\left\langle \left[ \theta _{1}\right] ,
\left[ \theta _{2}\right] ,\dots,\left[ \theta _{s}\right] \right\rangle =
{\rm Orb}\left\langle \left[ \vartheta _{1}\right] ,\left[ \vartheta
_{2}\right] ,\dots,\left[ \vartheta _{s}\right] \right\rangle .$$
\end{lemma}

Then, it exists a bijective correspondence between the set of ${\rm Aut}({\bf A})$-orbits on ${\rm T}_{s}\left( {\bf A}\right) $ and the set of
isomorphism classes of ${\rm E}\left( {\bf A},{{\bf V}}\right) $. Consequently we have a
procedure that allows us, given an assosymmetric algebra ${\bf A}'$ of
dimension $n-s$, to construct all non-split central extensions of ${\bf A}'$.

\; \;

{\centerline {\textsl{Procedure}}}

Let ${\bf A}'$ be an assosymmetric algebra of dimension $n-s $.

\begin{enumerate}
\item Determine ${\rm H^{2}}( {\bf A}',\mathbb {C}) $, ${\rm Ann}({\bf A}')$ and ${\rm Aut}({\bf A}')$.

\item Determine the set of ${\rm Aut}({\bf A}')$-orbits on ${\rm T}_{s}({\bf A}') $.

\item For each orbit, construct the assosymmetric algebra associated with a
representative of it.
\end{enumerate}

\subsection{Notations}
Let ${\bf A}$ be an assosymmetric algebra with
a basis $e_{1},e_{2},\dots,e_{n}$. Then by $\Delta _{ij}$\ we will denote the
assosymmetric bilinear form
$\Delta _{ij} \colon {\bf A}\times {\bf A}\longrightarrow \mathbb C$
with $\Delta _{ij}\left( e_{l},e_{m}\right) = \delta_{il}\delta_{jm}$.
Then the set $\left\{ \Delta_{ij}:1\leq i, j\leq n\right\} $ is a basis for the linear space of
the bilinear forms on ${\bf A}$. Then every $\theta \in
{\rm Z^2}\left( {\bf A} \right) $ can be uniquely written as $
\theta = \displaystyle \sum_{1\leq i,j\leq n} c_{ij}\Delta _{{i}{j}}$, where $
c_{ij}\in \mathbb C$.
Let us fix the following notations:
$$\begin{array}{lll}
{\mathcal A}^i_j& \mbox{---}& j\mbox{th }i\mbox{-dimensional one-generated nilpotent assosymmetric algebra.} \\
\end{array}$$

\subsection{The algebraic classification of  low dimensional one-generated nilpotent assosymmetric algebras}
In the present table (thanks to \cite{ikm19}) we have a description of all  $2$-, $3$- and $4$-dimensional  one-generated nilpotent assosymmetric algebras:

{\tiny
\begin{longtable}{l llllllllll}

${\mathcal A}^{2}_{01}$ &$:$ &  $e_1 e_1 = e_2$ &&\\ 
${\mathcal A}^{3}_{01}$ &$:$ & $e_1 e_1 = e_3$  & $e_2 e_1=e_3$ &  \\ 
${\mathcal A}^3_{02}(\alpha)$ &$:$& $e_1 e_1 = e_2$ & $e_1 e_2=e_3$ &  $e_2 e_1 = \alpha e_3$  \\  
${\mathcal A}^4_{01}$ &$:$& $e_1 e_1 = e_2$ & $e_1 e_2=e_4$ &  $e_2 e_1 = e_3$  \\  
${\mathcal A}^4_{02}$ &$:$& $e_1e_1=e_2$  &  $e_1e_2=e_4$  & $e_1e_3=e_4$ & $e_2e_1=e_3$ & $e_2e_2=-e_4$& $e_3e_1=-2e_4$\\  
${\mathcal A}^4_{03}$ &$:$& $e_1e_1=e_2$  &  $e_1e_3=e_4$   & $e_2e_1=e_3$ & $e_2e_2=-e_4$ & $e_3e_1=-2e_4$ \\ 
${\mathcal A}^4_{04}(\alpha)$ &$:$& $e_1e_1=e_2$ & $e_1e_2=e_3$ & $e_1e_3=(2-\alpha)e_4$   &  $e_2e_1=\alpha e_3$ & $e_2e_2=(\alpha^2-\alpha+1)e_4$ &  $e_3e_1=(2\alpha-1)e_4$ \\  
${\mathcal A}^4_{05}$ &$:$ & $e_1e_1=e_2$ & $e_1e_2=e_3$ & $e_1e_3=-2e_4$   & $e_2e_1=e_4$ & $e_2e_2=-e_4$ & $e_3e_1=e_4$ \\  
${\mathcal A}^4_{06}$ &$:$ & $e_1e_1=e_2$ & $e_1e_2=e_3$ & $e_1e_3=-3e_4$  & $e_2e_1=-e_3+e_4$ & $e_2e_2=-3e_4$ & $e_3e_1=3e_4$ \\ 

\end{longtable}
}

\begin{remark}
Note that, non-split central extension of a split algebra can not be a one-generated algebra. Hence, we will consider
central extensions only for non-split one-generated nilpotent algebras.
\end{remark} 
\section{Classification of $5$-dimensional one-generated nilpotent assosymmetric algebras}

\subsection{$2$-dimensional central extensions of $3$-dimensional one-generated algebras}
 The second cohomology spaces of algebras $ {\mathcal A}^3_{01}, {\mathcal A}^3_{02}(\alpha)$ given in \cite{ikm19}. Therefore, two dimensional central extensions of these algebras gives the following:

\begin{longtable}{{llllllll}}

${\mathcal A}^{5}_{01} $ &$:$ &$ e_1e_1=e_2  $&$ e_1e_2=e_4 $ &$ e_1e_3=e_5 $ \\&&$ e_2e_1=e_3 $  &$ e_2e_2=-e_5$   &$ e_3e_1=-2e_5$ &\\

${\mathcal A}^{5}_{02}(\alpha)$&$:$ &$ e_1e_1=e_2 $ &$ e_1e_2=e_3 $ &$ e_1e_3=(\alpha-2)e_5 $ \\&&$ e_2e_1=\alpha e_3 +e_4$  &$ e_2e_2=(\alpha-\alpha^2-1)e_5$   &$ e_3e_1=(1-2\alpha)e_5$
\end{longtable}

\subsection{Cohomology spaces  of $4$-dimensional one-generated assosymmetric algebras}

In the present table we collect all usefull information about ${\rm Z}^2, {\rm B}^2$ and ${\rm H}^2$ spaces for all $4$-dimensional one-generated algebras that were counted via code in \cite{km20}.

\begin{longtable}{|l c l |}

 \hline
 
${\rm Z^{2}}\left( {\mathcal A}^{4}_{01} \right)$ & $=$&  
$\left\langle \Delta_{11},\Delta_{12},\Delta_{21},\Delta_{13}+\Delta_{41}, 
\Delta_{14}-\Delta_{31}-\Delta_{41}, 
\Delta_{22}+2\Delta_{31}+\Delta_{41} \right\rangle$ \\

${\rm B^{2}}\left( {\mathcal A}^{4}_{01}\right)$ & $=$&  
$ \left\langle \Delta_{11},\Delta_{12},\Delta_{21}  \right\rangle$ \\

${\rm H^{2}}\left( {\mathcal A}^{4}_{01}\right)$ & $=$&  
$  \left\langle  [\Delta_{13}]+[\Delta_{41}],   [\Delta_{14}]-[\Delta_{31}]-[\Delta_{41}],   [\Delta_{22}]+2[\Delta_{31}]+[\Delta_{41}]  \right\rangle$ \\
  \hline

${\rm Z^{2}}\left( {\mathcal A}^{4}_{02} \right)$ & $=$&  
$ \left\langle  \Delta_{11},\Delta_{12},\Delta_{21},  \Delta_{13}-\Delta_{22}-2\Delta_{31}  \right\rangle$ \\ 
${\rm B^{2}}\left( {\mathcal A}^{4}_{02} \right)$ & $=$&  
$ \left\langle   \Delta_{11},\Delta_{21},\Delta_{12}+\Delta_{13}-\Delta_{22}-2\Delta_{31}  \right\rangle$ \\  
${\rm H^{2}}\left( {\mathcal A}^{4}_{02} \right)$ & $=$&  
$\left\langle  [\Delta_{12}] \right\rangle$  \\  
  \hline

${\rm Z^{2}}\left( {\mathcal A}^{4}_{03} \right)$ & $=$&  
$  \left\langle  \Delta_{11},\Delta_{12},\Delta_{21}, \Delta_{13}-\Delta_{22}-2\Delta_{31} \right\rangle$ \\

${\rm B^{2}}\left( {\mathcal A}^{4}_{03} \right)$ & $=$&  
$ \left\langle  \Delta_{11} ,\Delta_{21},\Delta_{13}-\Delta_{22}-2\Delta_{31}   \right\rangle$ \\

${\rm H^{2}}\left( {\mathcal A}^{4}_{03} \right)$ & $=$&  
$ \left\langle    [\Delta_{12}]  \right\rangle$ \\
\hline

${\rm Z^{2}}\left( {\mathcal A}^{4}_{04}(\alpha)_{\alpha \neq 1} \right)$ & $=$&  
$\left\langle  \Delta_{11},\Delta_{12}, \Delta_{21}, 
(2-\alpha)\Delta_{13}+ (\alpha^2-\alpha+1) \Delta_{22}+(2\alpha-1)\Delta_{31} \right\rangle$ \\

${\rm B^{2}}\left( {\mathcal A}^{4}_{04}(\alpha)_{\alpha \neq 1} \right)$ & $=$&  
$ \left\langle  \Delta_{11},\Delta_{12}+\alpha\Delta_{21},   (2-\alpha)\Delta_{13}+ (\alpha^2-\alpha+1)\Delta_{22}+(2\alpha-1)\Delta_{31}  \right\rangle$ \\

${\rm H^{2}}\left( {\mathcal A}^{4}_{04}(\alpha)_{\alpha \neq 1} \right)$ & $=$&  
$ \left\langle     [\Delta_{12}] \right\rangle$ \\
\hline

${\rm Z^{2}}\left( {\mathcal A}^{4}_{04}(1)  \right)$ & $=$&  
$ \left\langle   \Delta_{11},\Delta_{12}, \Delta_{21},
\Delta_{13}+\Delta_{22}+\Delta_{31}, 
\Delta_{14}+\Delta_{23}+\Delta_{32}+\Delta_{41} \right\rangle$ \\

${\rm B^{2}}\left( {\mathcal A}^{4}_{04}(1)  \right)$ & $=$&  $ \left\langle \Delta_{11},\Delta_{12}+\Delta_{21},  \Delta_{13}+\Delta_{22}+\Delta_{31}  \right\rangle$ \\

${\rm H^{2}}\left( {\mathcal A}^{4}_{04}(1)  \right)$ & $=$&  
$\left\langle    [\Delta_{21}], [\Delta_{14}]+[\Delta_{23}]+  [\Delta_{32}]+[\Delta_{41}] \right\rangle$ \\
\hline

${\rm Z^{2}}\left( {\mathcal A}^{4}_{05} \right)$ & $=$&  
$\left\langle   \Delta_{11}, \Delta_{12}, \Delta_{21},   
2\Delta_{13}+\Delta_{22}-\Delta_{31}   \right\rangle$ \\

${\rm B^{2}}\left( {\mathcal A}^{4}_{05} \right)$ & $=$&  
$ \left\langle  \Delta_{11},\Delta_{12},-2\Delta_{13}+\Delta_{21}-\Delta_{22}+\Delta_{31}         \right\rangle$ \\

${\rm H^{2}}\left( {\mathcal A}^{4}_{05} \right)$ & $=$&  
$ \left\langle    [\Delta_{21}] \right\rangle$\\
\hline

${\rm Z^{2}}\left( {\mathcal A}^{4}_{06} \right)$ & $=$&  
$ \left\langle   \Delta_{11}, \Delta_{12}, \Delta_{21},   \Delta_{13}+\Delta_{22}-\Delta_{31}   \right\rangle$ \\

${\rm B^{2}}\left( {\mathcal A}^{4}_{06} \right)$ & $=$&  
$\left\langle  \Delta_{11},\Delta_{12}-\Delta_{21},  -3\Delta_{13}+\Delta_{21}-3\Delta_{22}+3\Delta_{31}  \right\rangle$ \\

${\rm H^{2}}\left( {\mathcal A}^{4}_{06} \right)$ & $=$&  
$ \left\langle  [\Delta_{13}]+[\Delta_{22}]-[\Delta_{31}]  \right\rangle$ 
                     \\
                     \hline
\end{longtable}

\begin{remark}
Extensions of the algebras ${\mathcal A}^4_{02},$  ${\mathcal A}^4_{03},$ 
 ${\mathcal A}^4_{04}(\alpha)_{\alpha\neq 1},$  ${\mathcal A}^4_{05}$ and  ${\mathcal A}^4_{06}$ give algebras with $2$-dimensional annihilator. Then, in the following subsections we study the central extensions of the other algebras.
\end{remark}

  \subsection{Central extensions of ${\mathcal A}^{4}_{01}$}

Let us use the following notations: 
\[\nabla_1=[\Delta_{13}]+[\Delta_{41}], \ \nabla_2=[\Delta_{14}]-[\Delta_{31}]-[\Delta_{41}], \ \nabla_3=[\Delta_{22}]+2[\Delta_{31}]+[\Delta_{41}].
\]
The automorphism group of ${\mathcal A}^{4}_{01}$ consists of invertible matrices of the form

\[\phi=\left(
                             \begin{array}{cccc}
                               x & 0 & 0 & 0   \\
                               y & x^2 & 0 & 0  \\
                               z & xy & x^3 & 0 \\
                               t & xy & 0 & x^3                               \end{array}\right)
                               .\]

Since
\[
\phi^T
                           \left(\begin{array}{cccc}
                                0 & 0 & \alpha_1 & \alpha_2 \\
                                 0 & \alpha_3 & 0 & 0 \\
                                 -\alpha_2+2\alpha_3 & 0 & 0 & 0\\
                            \alpha_1-\alpha_2+\alpha_3 & 0 & 0 & 0\\
                             \end{array}
                           \right)\phi
                           =\left(\begin{array}{cccc}
                                \alpha^* & \alpha^{**} & \alpha^*_1 & \alpha^*_2 \\
                                 \alpha^{***} & \alpha^*_3 & 0 & 0 \\
                                 -\alpha^*_2+2\alpha^*_3 & 0 & 0 & 0\\
                            \alpha^*_1-\alpha^*_2+\alpha^*_3 & 0 & 0 & 0\\
                             \end{array}
                           \right),\]
we have that the action of ${\rm Aut} ({\mathcal A}^{4}_{01})$ on the subspace
$\langle  \sum\limits_{i=1}^3\alpha_i \nabla_i \rangle$
is given by
$\langle  \sum\limits_{i=1}^3\alpha^*_i \nabla_i \rangle,$ where

\[
\begin{array}{rclrclrcl}
\alpha^*_1&=&x^4 \alpha_1, &
\alpha^*_2&=&x^4 \alpha_2,& 
\alpha^*_3&=&x^4 \alpha_3.
\end{array}\]

\subsubsection{$1$-dimensional central extensions}
We have the following new cases:

\begin{enumerate}

\item If $\alpha_1\neq0,\alpha_2=0, \alpha_3=0,$ then $x=\frac{1}{\sqrt[4]{\alpha_1}},$  we have the representative $\langle \nabla_1\rangle;$

\item If $\alpha_2\neq0, \alpha_3=0,$ then $x=\frac{1}{\sqrt[4]{\alpha_2}}, \alpha=\frac{\alpha_1}{\alpha_2}$ we have the representative $\langle \alpha\nabla_1+\nabla_2\rangle;$ 

\item If $ \alpha_3\neq0,$ then $x=\frac{1}{\sqrt[4]{\alpha_3}}, \alpha=\frac{\alpha_1}{\alpha_3}, \beta=\frac{\alpha_2}{\alpha_3}$  we have the representative $\langle \alpha\nabla_1+\beta\nabla_2+\nabla_3\rangle.$

\end{enumerate}

From here, we have new $5$-dimensional one generated assosymmetric algebras constructed from  ${\mathcal
A}^{4}_{01}:$

\begin{longtable}{lllllll}

    ${\mathcal A}^{5}_{03} $&$:$&$ e_1e_1=e_2 $&$ e_1e_2=e_4 $&$ e_2e_1=e_3$&$ e_1e_3=e_5$&$ e_4e_1=e_5 $\\
    
   $ {\mathcal A}^{5}_{04}(\alpha) $&$:$&$ e_1e_1=e_2 $&$ e_1e_2=e_4 $&$ e_2e_1    =e_3$&$ e_1e_3=\alpha e_5 $\\ &&$e_1e_4=e_5 $&$ e_3e_1=-e_5 $&$ e_4e_1=(\alpha-1)e_5$\\
    
    ${\mathcal A}^{5}_{05}(\alpha,\beta) $&$:$&$ e_1e_1=e_2 $&$ e_1e_2=e_4 $&$ e_2e_1=e_3 $&$ e_1e_3=\alpha e_5 $ \\
    && $ e_1e_4=\beta e_5 $
    & $ e_2e_2=e_5 $&$  e_3e_1=(2-\beta)e_5 $& \multicolumn{2}{l}{$e_4e_1=(\alpha-\beta+1)e_5$}\\
\end{longtable}

\subsubsection{$2$-dimensional central extensions}
Consider the vector space generated by the following two cocycles

$$\begin{array}{rcl}
\theta_1 &=&\alpha_1\nabla_1+\alpha_2\nabla_2+\alpha_3\nabla_3,\\
\theta_2 &=&\beta_1\nabla_1+\beta_2\nabla_2.
\end{array}$$

Here we have the following cases:

\begin{enumerate}
\item If $\alpha_3=0,$ then we have the representative  $\langle \nabla_1,\nabla_2\rangle;$

\item If $\alpha_3\neq 0, \beta_1\neq0, \beta_2=0,$ then we have the representative $\langle \nabla_1, \alpha \nabla_2+\nabla_3\rangle;$

\item If $\alpha_3\neq 0, \beta_2 \neq 0,$ then we have the representative 
$\langle \alpha \nabla_1+\nabla_2, \beta \nabla_1   +\nabla_3\rangle.$

\end{enumerate}

We have the following new  $6$-dimensional one-generated nilpotent assosymmetric algebras constructed from  ${\mathcal
A}^{4}_{01}:$

\begin{longtable}{lllllll}

    ${\mathcal A}^6_{01}$&$:$&$ e_1e_1=e_2 $&$ e_1e_2=e_4$&$ e_1e_3=e_5$&$
    e_1e_4=e_6$\\
    &&$
    e_2e_1=e_3 $ &$e_3e_1=-e_6 $&$  e_4e_1=e_5-e_6$\\
     
    ${\mathcal A}^6_{02}(\alpha) $&$:$&$ e_1e_1=e_2 $&$ e_1e_2=e_4 $&$ e_1e_3=e_5 $&$ e_1e_4=\alpha e_6 $\\
    &&$  e_2e_1=e_3 $ &$ e_2e_2=e_6 $&$  e_3e_1=(2-\alpha) e_6$& \multicolumn{2}{l}{$ e_4e_1=e_5-(\alpha-1)e_6$}\\
     
    ${\mathcal A}^6_{03}(\alpha,\beta) $&$:$&$ e_1e_1=e_2 $&$ e_1e_2=e_4 $&$ e_1e_3=\alpha e_5+\beta e_6 $&$ e_1e_4=e_5 $\\
    &&$ e_2e_1=e_3$&$ e_2e_2=e_6 $&$ e_3e_1=- e_5+2e_6 $&  \multicolumn{2}{l}{$e_4e_1=(\alpha-1) e_5+(\beta+1)e_6$} \\
    
\end{longtable}

  \subsection{Central extensions of ${\mathcal A}^{4}_{04}(1)$}
Let us use the following notations: 

$$\nabla_1=[\Delta_{21}], \  \nabla_2=[\Delta_{14}]+[\Delta_{23}]+[\Delta_{32}]+[\Delta_{41}].$$

The automorphism group of ${\mathcal A}^{4}_{04}(1)$ consists of invertible matrices of the form

\[\phi= \left(\begin{array}{cccc}
                               x & 0   & 0  & 0  \\
                               y & x^2 & 0  & 0 \\
                               z & 2xy & x^3 & 0 \\
                               t & 2xz+y^2 &3yx^2 & x^4
                               \end{array}\right)
                               .\]

Since
\[
\phi^T
                           \left(\begin{array}{cccc}
                                0 & 0 & 0 & \alpha_2 \\
                                 \alpha_1 & 0 & \alpha_2 & 0 \\
                                 0 & \alpha_2 & 0 & 0\\
                            \alpha_2 & 0 & 0 & 0\\
                             \end{array}
                           \right)\phi
                           =\left(\begin{array}{cccc}
                                \alpha^{***} & \alpha^* & \alpha^{**} & \alpha^*_2 \\
                                 \alpha^*+\alpha^*_1 & \alpha^{**} & \alpha^*_2 & 0 \\
                                 \alpha^{**} & \alpha^*_2 & 0 & 0\\
                            \alpha^*_2 & 0 & 0 & 0\\
                             \end{array}
                           \right),\]
we have that the action of ${\rm Aut} ({\mathcal A}^{4}_{04}(1))$ on the subspace
$\langle  \sum\limits_{i=1}^2\alpha_i \nabla_i \rangle$
is given by
$\langle  \sum\limits_{i=1}^2\alpha^*_i \nabla_i \rangle,$ where

\[
\begin{array}{rclrcl}
\alpha^*_1&=&x^3 \alpha_1,  & 
\alpha^*_2&=&x^5 \alpha_2.
\end{array}\]

\subsubsection{$1$-dimensional central extensions}
Note that if $\alpha_2=0$ then we obtain algebras with 2-dimensional annihilator. Therefore, 
we have two representatives $\langle\nabla_2\rangle$ and $\langle\nabla_1+\nabla_2\rangle$ depending on whether $\alpha_1=0$ or not.

We have the following new $5$-dimensional nilpotent assosymmetric algebras constructed from  ${\mathcal
A}^{4}_{04}(1):$


\begin{longtable}{lllllll}
    ${\mathcal A}^{5}_{06}$ &$:$& $e_1e_1=e_2 $& $e_1e_2=e_3 $& $e_1e_3=e_4$ & $e_2e_1=e_3 $& $e_2e_2=e_4$ \\ 
    &  &$e_3e_1=e_4$ & $e_1e_4=e_5$ & $e_2e_3=e_5$ &    $ e_3e_2=e_5$ & $e_4e_1=e_5$\\
    
    ${\mathcal A}^{5}_{07} $&$ : $&$ e_1e_1=e_2 $&$ e_1e_2=e_3 $&$ e_1e_3=e_4 $&$ e_2e_1=e_3+e_5 $&$ e_2e_2=e_4 $\\ & & $ e_3e_1=e_4 $&$  e_1e_4=e_5 $&$ e_2e_3=e_5 $&$  e_3e_2=e_5 $&$ e_4e_1=e_5$\\

\end{longtable} 
 
 \subsubsection{$2$-dimensional central extensions}
 We have only one new $6$-dimensional nilpotent assosymmetric algebras constructed from  ${\mathcal
A}^{4}_{04}(1):$
 
\begin{longtable}{lllllll}
 ${\mathcal A}^6_{04}$ &$ : $&$ e_1e_1=e_2 $&$ e_1e_2=e_3 $&$ e_1e_3=e_4 $&$ e_1e_4=e_5 $&$ e_2e_1=e_3+e_6$ \\ 
    & &$ e_2e_2=e_4 $&$ e_2e_3=e_5 $&$ e_3e_1=e_4 $&$    e_3e_2=e_5 $&$ e_4e_1=e_5$

\end{longtable} 

 \subsection{Classification theorem}
 Summarizing results of the previous sections, we have the following theorem.

\begin{theoremA}
Let $\mathcal A$ be a $5$-dimensional complex one-generated nilpotent assosymmetric  algebra,
then $\mathcal A$ is isomorphic to an algebra  from the following list:
\end{theoremA}

\begin{longtable}{|l lllll|}

\hline

${\mathcal A}^{5}_{01} $ &$ e_1e_1=e_2  $&$ e_1e_2=e_4 $ &$ e_1e_3=e_5 $&& \\
&$ e_2e_1=e_3 $  &$ e_2e_2=-e_5$   &$ e_3e_1=-2e_5$ &&\\
\hline

${\mathcal A}^{5}_{02}(\alpha)  $&$ e_1e_1=e_2 $ &$ e_1e_2=e_3 $ &$ e_1e_3=(\alpha-2)e_5 $ &
 \multicolumn{2}{l|}{$ e_2e_1=\alpha e_3 +e_4$}  \\&  
 \multicolumn{2}{l}{$ e_2e_2=(\alpha-\alpha^2-1)e_5$}   &$ e_3e_1=(1-2\alpha)e_5$ &&\\

\hline 
${\mathcal A}^{5}_{03}$ & $e_1e_1=e_2 $&$ e_1e_2=e_4 $& $e_1e_3=e_5$ & $e_2e_1=e_3 $ &$ e_4e_1=e_5 $\\

\hline
${\mathcal A}^{5}_{04}(\alpha) $&$ e_1e_1=e_2 $&$ e_1e_2=e_4 $ &$  e_1e_3=\alpha e_5   $ &$ e_1e_4=e_5 $ & \\ 
&$ e_2e_1=e_3  $   &$ e_3e_1=-e_5  $&$ e_4e_1=(\alpha-1)e_5$ && \\

\hline
${\mathcal A}^{5}_{05}(\alpha,\beta)  $& $e_1e_1=e_2 $ &$ e_1e_2=e_4  $ &$ e_1e_3=\alpha e_5 $ &$ e_1e_4=\beta e_5  $& \\ 
& $  e_2e_1=e_3 $  &$ e_2e_2=e_5  $&$ e_3e_1=(2-\beta)e_5 $ & \multicolumn{2}{l|}{$  e_4e_1=(\alpha-\beta+1)e_5$}\\
\hline

 ${\mathcal A}^{5}_{06}  $&$ e_1e_1=e_2 $ &$ e_1e_2=e_3  $&$ e_1e_3=e_4  $&$ e_1e_4=e_5  $&$ e_2e_1=e_3  $\\ &$ e_2e_2=e_4
 $ &$ e_2e_3=e_5   $&$ e_3e_1=e_4  $  &$ e_3e_2=e_5$  &$ e_4e_1=e_5 $\\
\hline
    
 ${\mathcal A}^{5}_{07} $ &$  e_1e_1=e_2  $&$ e_1e_2=e_3  $&$ e_1e_3=e_4 $ &$ e_1e_4=e_5 $ &$ e_2e_1=e_3+e_5$ \\  &$ e_2e_2=e_4   $&$   e_2e_3=e_5   $&$  e_3e_1=e_4  $ &$  e_3e_2=e_5 $ &$ e_4e_1=e_5$  \\
\hline

\hline
 \end{longtable}

 \section{Classification of 6-dimensional one-generated nilpotent assosymmetric algebras}

\subsection{Cohomology spaces  of $5$-dimensional one-generated assosymmetric algebras}
All multiplication tables of $5$-dimensional one-generated nilpotent assosymmetric algebras is given in table in Theorem A (see, previous section).
All necessary information about 
coboundaries, cocycles and  second cohomology spaces of $5$-dimensional one-generated nilpotent assosymmetric algebras were calculated by the code in \cite{km20} and  given in the following table.

{\tiny
\begin{longtable}{|lll   |}

\hline
\multicolumn{3}{|c|}{\mbox{ Table B.  The list of cohomology spaces  of 5-dimensional one-generated assosymmetric algebras}}\\

\multicolumn{3}{|c|}{}\\

\hline 
		
${\rm Z^{2}}({\mathcal A}^{5}_{01})$& $=$ &$
		\Big\langle \begin{array}{l}\Delta_{11},\Delta_{12},\Delta_{21},\Delta_{13}+\Delta_{41},\Delta_{22}+2\Delta_{31}+\Delta_{41} 
			\Delta_{14}-\Delta_{31}-\Delta_{41} \end{array}\Big\rangle$\\
${\rm B^{2}}({\mathcal A}^{5}_{01})$& $=$ &$
			\Big\langle \begin{array}{l}\Delta_{11},\Delta_{12},\Delta_{21},
			\Delta_{13}-\Delta_{22}-2\Delta_{31} \end{array}\Big\rangle$\\
${\rm H^{2}}({\mathcal A}^{5}_{01})$& $=$ &$
			\Big\langle \begin{array}{l} [\Delta_{13}]+[\Delta_{41}],[\Delta_{14}]-[\Delta_{31}]+[\Delta_{41}]
			\end{array}\Big\rangle $\\
			\hline 

${\rm Z^{2}}(	{\mathcal A}^{5}_{02}(\alpha\neq1) )$& $=$ &$	         
            \Big\langle \begin{array}{l}\Delta_{11},\Delta_{12},\Delta_{21},\Delta_{13}+(1-\alpha)\Delta_{22}+(1-2\alpha)\Delta_{41},\\ \relax 
			\Delta_{14}-\Delta_{22}-2\Delta_{41},  
			\Delta_{22}+\Delta_{31}+(2-\alpha)\Delta_{41} \end{array}\Big\rangle $\\
			
${\rm B^{2}}(	{\mathcal A}^{5}_{02}(\alpha\neq1) )$& $=$ &$		 
			\Big\langle \begin{array}{l}\Delta_{11},\Delta_{12},\Delta_{21},
			(\alpha-2)\Delta_{13}+(\alpha-\alpha^2-1)\Delta_{22}+(1-2\alpha)\Delta_{31} \end{array}\Big\rangle$\\
			
${\rm H^{2}}(	{\mathcal A}^{5}_{02}(\alpha\neq1) )$& $=$ &$
			\Big\langle \begin{array}{l}[\Delta_{14}]-[\Delta_{22}]-2[\Delta_{41}], 			[\Delta_{22}]+[\Delta_{31}]+(2-\alpha)[\Delta_{41}]
			\end{array}\Big\rangle $\\
			\hline 
		
${\rm Z^{2}}(	{\mathcal A}^{5}_{02}(1) )$& $=$ &$	
            \Big\langle \begin{array}{l}\Delta_{11},\Delta_{12},\Delta_{21},\Delta_{13}-\Delta_{41},\Delta_{22}+\Delta_{31}+\Delta_{41}, 
			\Delta_{14}+\Delta_{31}-\Delta_{41},  
			\Delta_{15}-\Delta_{23}-\Delta_{32}+\Delta_{51} \end{array}\Big\rangle$\\
			
${\rm B^{2}}(	{\mathcal A}^{5}_{02}(1) )$& $=$ &$		
            \Big\langle \begin{array}{l}\Delta_{11},\Delta_{12},\Delta_{21},
			\Delta_{13}+\Delta_{22}+\Delta_{31} \end{array}\Big\rangle$\\
			
${\rm H^{2}}(	{\mathcal A}^{5}_{02}(1) )$& $=$ &$	
            \Big\langle \begin{array}{l} [\Delta_{13}]-[\Delta_{41}],[\Delta_{14}]+[\Delta_{31}]-[\Delta_{41}],			[\Delta_{15}]-[\Delta_{23}]-[\Delta_{32}]+[\Delta_{51}]
			\end{array}\Big\rangle$ \\
			\hline 
${\rm Z^{2}}({\mathcal A}^{5}_{03})$&$=$&$    
 \Big\langle \begin{array}{l}\Delta_{11},\Delta_{12},\Delta_{21},\Delta_{13}+\Delta_{41},  \relax 
			\Delta_{14}-\Delta_{31}-\Delta_{41},  \relax 
			\Delta_{22}+2\Delta_{31}+\Delta_{41} \end{array}\Big\rangle$  \\

${\rm B^{2}}({\mathcal A}^{5}_{03})$&$=$&$ 
\Big\langle \begin{array}{l}\Delta_{11},\Delta_{12},\Delta_{21},\Delta_{13}+\Delta_{41} \end{array}\Big\rangle$ \\
${\rm H^{2}}(		{\mathcal A}^{5}_{03})$&$=$&$   
\Big\langle \begin{array}{l} [\Delta_{14}]-[\Delta_{31}]-[\Delta_{41}],  [\Delta_{22}]+2[\Delta_{31}]+[\Delta_{41}] \end{array}\Big\rangle$ \\
			\hline

${\rm Z^{2}}(		{\mathcal A}^{5}_{04}(\alpha) )$&$=$&$
			\Big\langle \begin{array}{l}\Delta_{11},\Delta_{12},\Delta_{21},\Delta_{13}+\Delta_{14}-\Delta_{31}, \relax \Delta_{14}-\Delta_{31}-\Delta_{41},  \Delta_{14}+\Delta_{22}+\Delta_{31} \end{array}\Big\rangle $ \\
			
${\rm B^{2}}(	{\mathcal A}^{5}_{04}(\alpha) )  $&$=$&$ \Big\langle \begin{array}{l}\Delta_{11},\Delta_{12},\Delta_{21}, \relax
			\alpha\Delta_{13}+\Delta_{14}-\Delta_{31}+(\alpha-1)\Delta_{41} \end{array}\Big\rangle$ \\
${\rm H^{2}}(	{\mathcal A}^{5}_{04}(\alpha) ) $&$=$&$\Big\langle \begin{array}{l} [\Delta_{13}]+[\Delta_{41}], \alpha[\Delta_{13}]+2[\Delta_{14}]+[\Delta_{22}]+(\alpha-1)[\Delta_{41}] \end{array}\Big\rangle$\\
			\hline
			
${\rm Z^{2}}({\mathcal A}^{5}_{05}(\alpha,\beta))$ & $=$ &$
			\Big\langle \begin{array}{l}\Delta_{11},\Delta_{12},\Delta_{21},\Delta_{13}+\Delta_{14}-\Delta_{31}, \Delta_{14}-\Delta_{31}-\Delta_{41},  \Delta_{14}+\Delta_{22}+\Delta_{31} \end{array}\Big\rangle $\\ 
			
${\rm B^{2}}({\mathcal A}^{5}_{05}(\alpha,\beta)) $ &$=$&$		\Big\langle \begin{array}{l}\Delta_{11},\Delta_{12},\Delta_{21},
			\alpha\Delta_{13}+\beta\Delta_{14}+\Delta_{22}+(2-\beta)\Delta_{31}+(\alpha-\beta+1)\Delta_{41} \end{array}\Big\rangle$\\
${\rm H^{2}}({\mathcal A}^{5}_{05}(\alpha,\beta))$ &$=$&$		
			\Big\langle \begin{array}{l} [\Delta_{13}]+[\Delta_{14}]-[\Delta_{31}], [\Delta_{14}]-[\Delta_{31}]-[\Delta_{41}] \end{array}\Big\rangle $\\
\multicolumn{3}{|c|}{{\bf $ 	\alpha \neq \frac{1}{2} (\beta \pm\sqrt{-2+6 \beta -3 \beta ^2})$}}\\		\hline       
${\rm Z^{2}}({\mathcal A}^{5}_{05}( \alpha,\beta))$&$=$& $
  \Big\langle \begin{array}{l}\Delta _{11},\Delta _{12},\Delta _{21},\Delta _{14}-\Delta _{31}-\Delta _{41},\Delta _{13}+\Delta _{41},(2 \beta-1) \Delta _{15}+\\
  + (2\beta (\alpha-1) +1)\Delta _{23}+(\alpha+2 \beta^2-3 \beta+1)\Delta _{24}+(-2 \alpha \beta+3 \alpha+2 \beta^2-3 \beta+1)\Delta _{32}+\\
  +(-2 \alpha \beta+2 \alpha+2 \beta-1)\Delta _{42} +(2 \alpha-2 \beta+1)\Delta _{51},\Delta _{22}+2 \Delta _{31}+\Delta _{41} \end{array}\Big\rangle $\\
		&&\\
${\rm B^{2}}({\mathcal A}^{5}_{05}( \alpha,\beta))$&$=$& $
			\Big\langle \begin{array}{l}\Delta_{11},\Delta_{12},\Delta_{21},
		\alpha\Delta_{13}+\beta\Delta_{14}+\Delta_{22}
		+(2-\beta)\Delta_{31}+(\alpha-\beta+1)\Delta_{41} \end{array}\Big\rangle$\\
		&&\\
${\rm H^{2}}({\mathcal A}^{5}_{05}( \alpha,\beta))$&$=$& $
			\Big\langle \begin{array}{l} [\Delta _{14}]-[\Delta _{31}]-[\Delta _{41}],[\Delta _{13}]+[\Delta _{41}], 
			(2 \beta-1)[\Delta_{15}]+(2\beta (\alpha-1) +1)[\Delta_{23}] 
			+(\alpha+2 \beta^2-3 \beta+1) [\Delta_{24}]\\
			+(-2 \alpha \beta+3 \alpha+2 \beta^2-3 \beta+1)[\Delta_{32}]
		+(-2 \alpha \beta+2 \alpha+2 \beta-1) [\Delta_{42}] 
		+(2\alpha-2\beta +1) [\Delta_{51}] \end{array}\Big\rangle$ \\
	\multicolumn{3}{|c|}{$	\alpha=\frac{1}{2}(\beta \pm \sqrt{-2+6 \beta -3 \beta ^2})$ and 
	$(\alpha,  \beta)\neq (0,\frac{1}{2})$} \\		\hline
			${\rm Z^{2}}(	{\mathcal A}^{5}_{05}(0,\frac{1}{2})) $& $=$& 
			$\Big\langle \begin{array}{l}\Delta_{11},\Delta_{12},\Delta_{21},\Delta_{13}+\Delta_{41}, \Delta_{14}-\Delta_{31}-\Delta_{41}, \\ \relax 2\Delta_{15}-3\Delta_{23}-2\Delta_{24}-3\Delta_{32}+\Delta_{42}-4\Delta_{51}, \Delta_{22}+2\Delta_{31}+\Delta_{41} \end{array}\Big\rangle$\\
		${\rm B^{2}}(	{\mathcal A}^{5}_{05}(0,\frac{1}{2})) $	& 
			$=$&$\Big\langle \begin{array}{l}\Delta_{11},\Delta_{12},\Delta_{21},			\Delta_{14}+2\Delta_{22}+3\Delta_{31}+\Delta_{41}\end{array}\Big\rangle$\\
		${\rm H^{2}}(	{\mathcal A}^{5}_{05}(0,\frac{1}{2})) $	&$=$&$
			\Big\langle \begin{array}{l} [\Delta_{14}]-[\Delta_{31}]-[\Delta_{41}], [\Delta_{13}]+[\Delta_{41}],
			2[\Delta_{15}]-3[\Delta_{23}]-2[\Delta_{24}]-3[\Delta_{32}]+[\Delta_{42}]-4[\Delta_{51}]
			\end{array}\Big\rangle$ \\
			\hline  
			

${\rm Z^{2}}(	{\mathcal A}^{5}_{06}))$ &$=$				
	 &  $
			\Big\langle \begin{array}{l}\Delta_{11},\Delta_{12},\Delta_{21},\Delta_{13}+\Delta_{22}+\Delta_{31},		\Delta_{14}+\Delta_{23}+\Delta_{32}+\Delta_{41}, \\ \relax 
			\Delta_{15}+\Delta_{24}+\Delta_{33}+\Delta_{42}+\Delta_{51} \end{array}\Big\rangle$ \\
			
${\rm B^{2}}(	{\mathcal A}^{5}_{06}))$ &$=$ &  $ 
			\Big\langle \begin{array}{l}\Delta_{11},\Delta_{12}+\Delta_{21},		\Delta_{13}+\Delta_{22}+\Delta_{31},	\Delta_{14}+\Delta_{23}+\Delta_{32}+\Delta_{41} \end{array}\Big\rangle$\\

${\rm H^{2}}(	{\mathcal A}^{5}_{06}))$ &$=$&  $
			\Big\langle \begin{array}{l}[\Delta_{21}],[\Delta_{15}]+[\Delta_{24}]+[\Delta_{33}]+[\Delta_{42}]+[\Delta_{51}]  \end{array}\Big\rangle$ \\
			\hline 
			
${\rm Z^{2}}({\mathcal A}^{5}_{07})$& $=$ &$
			\Big\langle \begin{array}{l}\Delta_{11},\Delta_{12},\Delta_{21},\Delta_{13}+\Delta_{22}+\Delta_{31},
			\Delta_{14}+\Delta_{23}+\Delta_{32}+\Delta_{41}, \\ \relax 
			\Delta_{15}+2\Delta_{22}+\Delta_{24}+3\Delta_{31}+\Delta_{33}+\Delta_{42}+\Delta_{51} \end{array}\Big\rangle $\\
			
${\rm B^{2}}({\mathcal A}^{5}_{07})$& $=$ &$ 
			\Big\langle \begin{array}{l}\Delta_{11},\Delta_{12}+\Delta_{21},			\Delta_{13}+\Delta_{22}+\Delta_{31},\Delta_{14}+\Delta_{21}+\Delta_{23}+\Delta_{32}+\Delta_{41} \end{array}\Big\rangle$\\
	
${\rm H^{2}}({\mathcal A}^{5}_{07})$& $=$ &$
			\Big\langle \begin{array}{l}[\Delta_{21}], 
			[\Delta_{15}]+2[\Delta_{22}]+[\Delta_{24}]+3[\Delta_{31}]+[\Delta_{33}]+[\Delta_{42}]+[\Delta_{51}]
			\end{array}\Big\rangle $\\
			\hline

\end{longtable}}

\begin{remark}
Extensions of the algebras ${\mathcal A}^5_{01},$ ${\mathcal A}^5_{02}(\alpha)_{\alpha\neq 1},$ ${\mathcal A}^5_{03},$ ${\mathcal A}^5_{04}(\alpha)$ and 
${\mathcal A}^5_{05}(\alpha, \beta)_{\alpha \neq \frac{1}{2} (\beta \pm\sqrt{-2+6 \beta -3 \beta ^2})}$ 
give algebras with $2$-dimensional annihilator. Then, in the following subsections we study the central extensions of the other algebras.
\end{remark}

  \subsection{Central extensions of ${\mathcal A}^{5}_{02}(1)$}

Let us use the following notations: 
\[\nabla_1=[\Delta_{13}]-[\Delta_{41}], \ \nabla_2=[\Delta_{14}]+[\Delta_{31}]-[\Delta_{41}],
\nabla_3=[\Delta_{15}]-[\Delta_{23}]-[\Delta_{32}]+[\Delta_{51}].
\]
The automorphism group of ${\mathcal A}^{5}_{07}$ consists of invertible matrices of the form

\[\phi= \left(\begin{array}{ccccc}
                               x & 0   & 0  & 0 & 0 \\
                               y & x^2 & 0 & 0 & 0 \\
                               z & 2xy &x^3 & 0 & 0 \\
                               t & xy & 0 & x^3 & 0\\
                               w & -y^2-2xz & -3x^2y  & 0 & x^4 
                               \end{array}\right)
                               .\]

Since
\[
\phi^T
                           \left(\begin{array}{ccccc}
                                0 & 0 & \alpha_1 & \alpha_2 & \alpha_3 \\
                                 0&  0 & -\alpha_3 & 0 & 0 \\
                                 \alpha_2 & -\alpha_3 & 0 &  0 & 0\\
                             -\alpha_1-\alpha_2 & 0 & 0 & 0 & 0\\
                            \alpha_3 & 0 & 0 & 0 & 0\\
                             \end{array}
                           \right)\phi
                           = \left(\begin{array}{ccccc}
                                \alpha^{****} & \alpha^{***} & \alpha^*_1+\alpha^* & \alpha^*_2 & \alpha^*_3 \\
                                 \alpha^{**}&  \alpha^* & -\alpha^*_3 & 0 & 0 \\
                                 \alpha^*_2+\alpha^{*} & -\alpha^*_3 & 0 &  0 & 0\\
                             -\alpha^*_1-\alpha^*_2 & 0 & 0 & 0 & 0\\
                            \alpha^*_3 & 0 & 0 & 0 & 0\\
                             \end{array}
                           \right)\]
we have that the action of ${\rm Aut} ({\mathcal A}^{5}_{07})$ on the subspace
$\langle  \sum\limits_{i=1}^3\alpha_i \nabla_i \rangle$
is given by
$\langle  \sum\limits_{i=1}^3\alpha^*_i \nabla_i \rangle,$ where

\[
\begin{array}{rclrclrcl}
\alpha^*_1&=&x^4\alpha_1, &
\alpha^*_2&=& x^4\alpha_2, &
\alpha^*_3&=& x^5\alpha_3.

\end{array}\]

We have the following case:

\begin{enumerate}

\item If $\alpha_2\neq0,$ then choosing $x=\frac{\alpha_2}{\alpha_3}$   we have the representative $\langle\alpha\nabla_1+\nabla_2+\nabla_3\rangle;$

\item If $\alpha_2=0,$ we have two representatives $\langle\nabla_3\rangle$ and $\langle\nabla_1+\nabla_3\rangle$ depending on whether $\alpha_1=0$ or not.
\end{enumerate}

Consequently, we have the following algebras from ${\mathcal A}^{5}_{02}(1):$

{\tiny
\begin{longtable}{lllllllllllllllllllllllll}

  $ {\mathcal A}^{6}_{05}(\alpha)$ &$: $&$  
  e_1e_1=e_2  $&$ e_1e_2=e_3   $&$  e_1e_3=-e_5+\alpha e_6  $&$ e_1e_4= e_6  $ &$ e_1e_5=e_6$&$ e_2e_1=e_3 +e_4$ \\
   &&$ e_2e_2=-e_5   $&$ e_2e_3=-e_6 $&$ e_3e_1=-e_5+ e_6  $
       &$ e_3e_2=-e_6$& $ e_4e_1=-(\alpha+1)  e_6  $&$ e_5e_1=e_6 $\\
   
${\mathcal A}^{6}_{06}$  &$: $&
$e_1e_1=e_2  $&$ e_1e_2=e_3   $&$  e_1e_3=-e_5+ e_6  $&$ e_1e_5=e_6  $ &$ e_2e_1=e_3 +e_4$  &$ e_2e_2=-e_5  $&\\&&$ e_2e_3=-e_6 $&$ e_3e_1=-e_5  $&$ e_3e_2=-e_6  $&$ e_4e_1=- e_6  $&$ e_5e_1=e_6 $\\
   
$ {\mathcal A}^{6}_{07}$ &$: $&$  e_1e_1=e_2  $&$ e_1e_2=e_3   $&$  e_1e_3=-e_5  $&$ e_1e_5=e_6  $  &$ e_2e_1=e_3 +e_4$ 
& $ e_2e_2=-e_5  $ \\&&$ e_2e_3=-e_6 $&$ e_3e_1=-e_5  $&$ e_3e_2=-e_6  $&$ e_5e_1=e_6 $\\
\end{longtable}
}

\subsection{Central extensions of ${\mathcal A}^{5}_{05}(  \alpha,\beta)$}  
Here we will consider the special cases for 
$\alpha=\frac{1}{2}(\beta \pm \sqrt{-2+6\beta -3 \beta ^2}).$

	The automorphism group of ${\mathcal A}^{5}_{05}(  \alpha,\beta)$ consists of invertible matrices of the form	
	$$\phi=\left(
\begin{array}{ccccc}
 x & 0 & 0 & 0 & 0 \\
 \frac{y}{x} & x^2 & 0 & 0 & 0 \\
 z & y & x^3 & 0 & 0 \\
 t & y & 0 & x^3 & 0 \\
 v & \frac{ x^3((2-\beta+\alpha)z   +(1+\alpha)t)+y^2}{x^2} & (\alpha-2 \beta+4) x y & (\alpha+\beta+1) x y & x^4 \\
\end{array}
\right).$$

Let use the following notations: 
\[\nabla_1=	[\Delta_{14}]-[\Delta_{31}]-[\Delta_{41}],
 \nabla_2=[\Delta_{13}]+[\Delta_{41}],\]
\[\nabla_3=	(2 \beta-1)[\Delta_{15}]+(2\alpha\beta -2\beta +1)[\Delta_{23}] +(\alpha+2 \beta^2-3 \beta+1) [\Delta_{24}]+\] 
\[ (3 \alpha-2 \alpha \beta+2 \beta^2-3 \beta+1)[\Delta_{32}]
		+(2 \alpha-2 \alpha \beta+2 \beta-1) [\Delta_{42}] 
		+(2\alpha-2\beta +1) [\Delta_{51}]. \]

So, {\tiny$$\phi^T\left(
\begin{array}{ccccc}
 0 & 0 & \alpha _2 & \alpha _1 &(2\beta-1)\alpha _3 \\
 0 & 0 & (2\alpha\beta -2\beta +1)\alpha _3 & (\alpha+2 \beta^2-3 \beta+1)\alpha _3 & 0 \\
 -\alpha _1 &(3 \alpha-2 \alpha \beta+2 \beta^2-3 \beta+1)\alpha _3 & 0 & 0 & 0 \\
 -\alpha_1+\alpha_2 & (2 \alpha-2 \alpha \beta+2 \beta-1)\alpha _3& 0 & 0 & 0 \\
(2\alpha-2\beta +1) \alpha _3 & 0 & 0 & 0 & 0 \\
\end{array}
\right)\phi=$$
$$=\left(
\begin{array}{ccccc}
 \alpha^{****} & \alpha^{***} &\alpha\alpha^*+\alpha^* _2 &\beta\alpha^*+ \alpha^* _1 &(2\beta-1)\alpha^* _3 \\
 \alpha^{**} &\alpha^* & (2\alpha\beta -2\beta +1) \alpha^* _3 &  (\alpha+2 \beta^2-3 \beta+1)\alpha^* _3 & 0 \\
 (2-\beta)\alpha^*-\alpha^* _1& (3 \alpha-2 \alpha \beta+2 \beta^2-3 \beta+1)\alpha^* _3 & 0 & 0 & 0 \\
 (1-\beta+\alpha)\alpha^*-\alpha^*_1++\alpha^*_2 &  (2 \alpha-2 \alpha \beta+2 \beta-1)\alpha^* _3& 0 & 0 & 0 \\
 (2\alpha-2\beta +1) \alpha^* _3 & 0 & 0 & 0 & 0 \\
\end{array}
\right)$$}
	
we have that the action of ${\rm Aut} ({\mathcal A}^{5}_{05}(\alpha,\beta))$ on the subspace
$\langle  \sum\limits_{i=1}^3\alpha_i \nabla_i \rangle$
is given by
$\langle  \sum\limits_{i=1}^3\alpha^*_i \nabla_i \rangle,$
where

$$\begin{array}{lclcl}
   
 \alpha^* _1&=&x^4\alpha _1-\beta(\beta-2)  (4\beta-2\alpha-2)\alpha _3 x^2 y, \\
 \alpha _2^*&=&x^4\alpha_2 -\left(\beta (\beta -2) (2 \beta -1)+\alpha   (2 \beta ^2-4 \beta +3 )\right)\alpha _3 x^2 y, \\ 
 \alpha _3^*&=&x^5\alpha _3. \\
 \end{array}$$

We are interested only in the cases with $\alpha_3 \neq 0.$ Now we obtain the following cases:

\begin{enumerate}

\item For $\beta (\beta -2) (2 \beta -1)+\alpha  \left(2 \beta ^2-4 \beta +3\right)\neq 0$ :

\begin{enumerate}
    \item If 
$2 \beta  (\beta -2)   (   2 \beta-\alpha -1)\alpha_2=\alpha_1 \left(\beta (\beta -2) (2 \beta -1)+\alpha   (2 \beta ^2-4 \beta +3 )\right),$ then by choosing $x=\frac{1}{\sqrt[5]{\alpha_3}}$ and $y=\frac{\alpha_2 x^2}{ \beta (\beta -2) (2 \beta -1)+\alpha  \left(2 \beta ^2-4 \beta +3\right) },$ we have the representative $\langle\nabla_3\rangle;$

\item If $2 \beta (\beta -2)    ( 2 \beta-\alpha  -1)\alpha_2\neq\alpha_1 \left(\beta (\beta -2) (2 \beta -1)+\alpha  (2 \beta ^2-4 \beta +3 )\right),$ then by choosing 
$x=\frac{ \alpha_1\left(\alpha  (2 \beta ^2-4 \beta +3 )+\beta(\beta -2)    (2 \beta -1)\right)+2 \alpha_2 \beta(\beta -2)    (\alpha -2 \beta +1)}{ \beta (\beta -2) (2 \beta -1)+\alpha   (2 \beta ^2-4 \beta +3 ) }$ and $y=\frac{\alpha_2 x^2}{ \beta (\beta -2) (2 \beta -1)+\alpha  \left(2 \beta ^2-4 \beta +3\right) },$  and we have the representative $\langle\nabla_1+\nabla_3\rangle.$

\end{enumerate}
From the above cases we have new parametric algebras:


 {\tiny\begin{longtable}{llllllllll} 
    ${\mathcal A}^6_{i}(\beta) $&$:$& 
    $e_1 e_1=e_2 $&$  e_1 e_2=e_4$ & $ e_1 e_3=\alpha e_5$ \\
    &&$  e_1 e_4=\beta e_5 $ &$ e_1 e_5=(2 \beta-1) e_6 $& $ e_2 e_1=e_3$ \\
    &&$  e_2 e_2=e_5 $&$  e_2 e_3= (2 \alpha \beta-2 \beta+1)e_6  $ &$  e_2 e_4= (\alpha+2 \beta^2-3 \beta+1)e_6 $\\ 
    &&$ e_3 e_1=(2-\beta) e_5 $&$  e_3 e_2=(3 \alpha-2 \alpha \beta+2 \beta^2-3 \beta+1)e_6 $&$   e_4 e_1= (\alpha-\beta+1)e_5$
    \\&&$ e_4 e_2=( 2 \alpha-2 \alpha \beta+2 \beta-1)e_6$  &$  e_5 e_1= (2 \alpha-2 \beta+1)e_6$ &&\\ 
    
     ${\mathcal A}^6_{i+1}(\beta) $&$:$& 
    $e_1 e_1=e_2 $&$  e_1 e_2=e_4$ & $ e_1 e_3=\alpha e_5$ \\
    &&$  e_1 e_4=\beta e_5+e_6 $ &$ e_1 e_5=(2 \beta-1) e_6 $& $ e_2 e_1=e_3$ \\
    &&$  e_2 e_2=e_5 $&$  e_2 e_3= (2 \alpha \beta-2 \beta+1)e_6  $ &$  e_2 e_4= (\alpha+2 \beta^2-3 \beta+1)e_6 $\\ 
    &&$ e_3 e_1=(2-\beta) e_5-e_6 $&$  e_3 e_2=( 3 \alpha-2 \alpha \beta+2 \beta^2-3 \beta+1)e_6 $&$   e_4 e_1= (\alpha-\beta+1)e_5-e_6$
    \\&&$ e_4 e_2=( 2 \alpha-2 \alpha \beta+2 \beta-1)e_6$  &$  e_5 e_1= (2 \alpha-2 \beta+1)e_6$ &&
\end{longtable} }
where $i=08$ for $\alpha=\frac{1}{2}(\beta +\sqrt{-2+6\beta -3 \beta ^2})$ with ${\beta \not \in \{1,\frac{3}{2}\}}$
and 
where $i=10$ for $\alpha=\frac{1}{2}(\beta-\sqrt{-2+6\beta -3 \beta ^2})$ with ${\beta \neq \frac{1}{2}}.$

\item The condition  $\beta=1,$  for $\alpha=\frac{1}{2}(\beta + \sqrt{-2+6\beta -3 \beta ^2})$ gives  $ \alpha=1,$ that is ${\mathcal A}^{5}_{05}(1,1)$. 
The base of the second cohomology of this algebra spanned by elements: 
\[     \nabla_1=[\Delta_{14}]-[\Delta_{31}]-[\Delta_{41}],\
    \nabla_2=[\Delta_{13}]+[\Delta_{41}],\
    \nabla_3=[\Delta_{15}]+[\Delta_{23}]+[\Delta_{24}]+[\Delta_{32}]+[\Delta_{42}]+[\Delta_{51}]. \] 
    Since
\[
                \phi^T \left(
\begin{array}{ccccc}
 0 & 0 & \alpha _2 & \alpha _1 & \alpha _3 \\
 0 & 0 & \alpha _3 & \alpha _3 & 0 \\
 -\alpha_1 & \alpha _3 & 0 & 0 & 0 \\
 \alpha _2-\alpha _1 & \alpha _3 & 0 & 0 & 0 \\
 \alpha _3 & 0 & 0 & 0 & 0 \\
\end{array}
\right)\phi=
            \left(
\begin{array}{ccccc}
\alpha^{****} & \alpha^{**} & \alpha^{*}+\alpha^* _2 &\alpha^{*}+\alpha^* _1 & \alpha^* _3 \\
\alpha^{***} & \alpha^{*} & \alpha^* _3 & \alpha^* _3 & 0 \\
 \alpha^{*}-\alpha^* _1 & \alpha^* _3 & 0 & 0 & 0 \\
 \alpha^{*}-\alpha^* _1+\alpha^* _2 & \alpha^* _3 & 0 & 0 & 0 \\
 \alpha^* _3 & 0 & 0 & 0 & 0 \\
\end{array}
\right)        \]
we have that the action of ${\rm Aut} ({\mathcal A}^{5}_{05}(1,1))$ on the subspace
$\langle  \sum\limits_{i=1}^3\alpha_i \nabla_i \rangle$
is given by
$\langle  \sum\limits_{i=1}^3\alpha^*_i \nabla_i \rangle,$ where
$$
\begin{array}{ccc}
 \alpha _1^*=x^4\alpha _1, &
 \alpha _2^*=x^4\alpha _2, &
 \alpha _3^*=\alpha _3 x^5. \\
\end{array}
$$

We are interested only in $\alpha_3\neq0,$ then we have the following cases:
\begin{enumerate}

\item If $\alpha_2\neq0,$ then for $x=\frac{\alpha_2}{\alpha_3},$  $ \alpha=\frac{\alpha_1}{\alpha_2}$  we have the representative $\langle\alpha\nabla_1+\nabla_2+\nabla_3\rangle.$ 

\item If $\alpha_2=0,$ then also we have two cases:
\begin{enumerate} 
\item  If $\alpha_1\neq0,$ then
$x=\frac{\alpha_1}{\alpha_3}, $ and we have the representative $\langle \nabla_1+\nabla_3\rangle;$
\item  If $\alpha_1=0,$ then
$x=\frac{1}{\sqrt[5]{\alpha_3}}, $ and we have the representative $\langle\nabla_3\rangle;$
\end{enumerate}
\end{enumerate}

Consequently, we have the following algebras from ${\mathcal A}^{5}_{05}(1,1):$
${\mathcal A}^6_{08}(1),$ ${\mathcal A}^6_{09}(1)$ and 

\begin{longtable}{llllllll}
      ${\mathcal A}^6_{12}(\alpha)$&$:$ & $e_1e_1=e_2$ & $e_1e_2=e_4$  & $e_1e_3= e_5+ e_6$ & $e_1e_4= e_5+\alpha e_6 $ &  $e_1e_5=e_6$ \\&&$e_2e_1=e_3$ & $e_2e_2=e_5 $& $e_2e_3=e_6$ & $e_2e_4=e_6$ & $e_3e_1=e_5-\alpha e_6$ \\&&  $e_3e_2=e_6$ & $e_4e_1=e_5+(1-\alpha) e_6$  & $e_4e_2=e_6$ & $e_5e_1=e_6$ &\\
      
      
\end{longtable}

\item The condition $\beta=\frac{3}{2}$ gives $\alpha=1$ for $\alpha=\frac{1}{2}(\beta + \sqrt{-2+6\beta -3 \beta ^2}),$  that is ${\mathcal A}^{5}_{05}(1,\frac{3}{2}).$ 
 So, the second cohomology space of ${\mathcal A}^{5}_{05}(1,\frac{3}{2})$ spanned by elements: 
\[     \nabla_1=[\Delta_{14}]-[\Delta_{31}]-[\Delta_{41}],\
    \nabla_2=[\Delta_{13}]+[\Delta_{41}],\
    \nabla_3=2[\Delta_{15}]+[\Delta_{23}]+2[\Delta_{24}]+[\Delta_{32}]+[\Delta_{42}]. \]



    Since
\[
\phi^T
                           \left(
\begin{array}{ccccc}
 0 & 0 & \alpha _2 &  \alpha _1 & 2 \alpha _3 \\
 0 & 0 & \alpha _3 & 2 \alpha _3 & 0 \\
 -\alpha _1 &  \alpha _3 & 0 & 0 & 0 \\
 \alpha _2-\alpha _1 &  \alpha _3 & 0 & 0 & 0 \\
 0 & 0 & 0 & 0 & 0 \\
\end{array}
\right)\phi
                           =  \left(
\begin{array}{ccccc}
 \alpha^{****} & \alpha^{***} &\alpha^* _2+\alpha^{*} & \alpha^*_1+3  \alpha^{*} & 2\alpha^* _3 \\
  \alpha^{**} &2 \alpha^{*} & \alpha^* _3 & 2\alpha^* _3 & 0 \\
 \alpha^{*}-\alpha^* _1 & \alpha^* _3 & 0 & 0 & 0 \\
 
    \alpha^*_2-\alpha^*_1  +\alpha^* & \alpha^*_3 & 0 & 0 & 0 \\

 0 & 0 & 0 & 0 & 0 \\
\end{array}
\right)  \]

we have that the action of ${\rm Aut} ({\mathcal A}^{5}_{05}(1,\frac{3}{2}))$ on the subspace
$\langle  \sum\limits_{i=1}^3\alpha_i \nabla_i \rangle$
is given by
$\langle  \sum\limits_{i=1}^3\alpha^*_i \nabla_i \rangle,$ where

\[
\begin{array}{rclrclrcl}
\alpha^*_1&=&x^4 \alpha_1+\frac{3}{2} x^3 y \alpha_3, & 
\alpha^*_2&=&x^4 \alpha_2, &
\alpha^*_3&=& x^5 \alpha_3.\\
\end{array}\]
 
 Since $\alpha_3\neq0,$ and  choosing $y=-\frac{2x^2\alpha_1}{3\alpha_3},$ we have the representatives $\langle\nabla_3\rangle$ and $\langle\nabla_2+\nabla_3\rangle,$ depending on whether $\alpha_2=0$ or not.
 
  We have the following new $6$-dimensional  algebras constructed from ${\mathcal A}^{5}_{05}(1,\frac{3}{2}):$ 
  ${\mathcal A}^6_{08}(\frac{3}{2}) $ and  
  
  \begin{longtable}{llllllllllll}
   ${\mathcal A}^6_{13}:$ &$ e_1e_1=e_2  $&$ e_1e_2=e_4  $&$ e_1e_3= e_5+e_6 $&$ e_1e_4=\frac{3}{2}e_5$&$ e_1e_5=2 e_6$\\
     &$ e_2e_1=e_3$ &$ e_2e_2=e_5  $&$e_2e_3=e_6  $&$ e_2e_4=2e_6  $   &$  e_3e_1=\frac{1}{2}e_5  $\\      &$ e_3e_2=e_6 $ &$ e_4e_1=e_6  $ &$ e_4e_2=e_6  $        &  \\     \\
       
    \end{longtable}

\end{enumerate}

\subsection{Central extensions of ${\mathcal A}^{5}_{05}(0,\frac{1}{2})$} 
 If $\beta=\frac{1}{2}$ for $\alpha=\frac{1}{2}(\beta - \sqrt{-2+6\beta -3 \beta ^2})$ gives $ \alpha=0,$ that is ${\mathcal A}^{5}_{05}(0,\frac{1}{2})$. 
 So, the second cohomology space of ${\mathcal A}^{5}_{05}(0,\frac{1}{2})$ spanned by elements: 
\[\nabla_1=[\Delta_{14}]-[\Delta_{31}]-[\Delta_{41}], 
\nabla_2=[\Delta_{13}]+[\Delta_{41}],
\nabla_3=2[\Delta_{15}]-3[\Delta_{23}]-2[\Delta_{24}]-3[\Delta_{32}]+[\Delta_{42}]-4[\Delta_{51}].\]

    Since
\[
                \phi^T \left(
\begin{array}{ccccc}
 0 & 0 & \alpha _2 & \alpha _1 & 2 \alpha _3 \\
 0 & 0 & -3 \alpha _3 & -2 \alpha _3 & 0 \\
 -\alpha _1 & -3 \alpha _3 & 0 & 0 & 0 \\
 \alpha _2-\alpha _1 & \alpha _3 & 0 & 0 & 0 \\
 -4 \alpha _3 & 0 & 0 & 0 & 0 \\
\end{array}
\right)\phi
                           = \left(\begin{array}{ccccc}
                                \alpha^{****} &  \alpha^{**} & \alpha^*_2 & \alpha^*_1+\alpha^* & 2\alpha^*_3 \\
                                 \alpha^{***} & 2\alpha^{*} &   -3\alpha^*_3 & -2\alpha^*_3 & 0 \\
                                 -\alpha^*_1+3\alpha^* &  -3\alpha^*_3 & 0 &  0 & 0\\
                            \alpha^*_2-\alpha^*_1 +\alpha^*&  \alpha^*_3 & 0 & 0 & 0\\
                            -4\alpha^*_3 & 0 & 0 & 0 & 0\\
                             \end{array}
                           \right)\]
we have that the action of ${\rm Aut} ({\mathcal A}^{5}_{05}(0,\frac{1}{2}))$ on the subspace
$\langle  \sum\limits_{i=1}^3\alpha_i \nabla_i \rangle$
is given by
$\langle  \sum\limits_{i=1}^3\alpha^*_i \nabla_i \rangle,$ where
\[
\begin{array}{rclrclrcl}

\alpha^*_1&=&x^4 \alpha_2+ \frac{9}{2} x^3 y\alpha_3, & \
\alpha^*_2&=&x^4 \alpha_1+ 3 x^3 y\alpha_3, & \
\alpha^*_3&=& x^5 \alpha_3.\\
\end{array}\]

We are interested in $\alpha_3\neq0,$ then we have the following cases:

\begin{enumerate}

\item If $3\alpha_1-2\alpha_2=0,$ then $x=\frac{1}{\sqrt[5]{\alpha_3}}$ and $ y=-\frac{x\alpha_1 }{3\alpha_3},$  we have the representative $\langle\nabla_3\rangle;$

\item If $3\alpha_1-2\alpha_2\neq0,$ then $x=\frac{-3 \alpha_1+2 \alpha_2}{{2\alpha_3}}, y=-\frac{x\alpha_1}{3\alpha_3}$ and we have the representative $\langle \nabla_2+\nabla_3\rangle.$

\end{enumerate}

We have the following new $6$-dimensional  algebras constructed from ${\mathcal A}^{5}_{05}(0,\frac{1}{2}):$ 

\begin{longtable}{lllllllllllll}     
   ${\mathcal A}^6_{14}:  $&$ e_1e_1=e_2  $&$ e_1e_2=e_4  $&$ e_1e_4=\frac{1}{2}e_5 $
     
     &$ e_1e_5=2e_6$&$ e_2e_1=e_3  $\\&$ e_2e_2=e_5  $&$ e_2e_3=-3e_6  $ 
     &$ e_2e_4=-2e_6  $&$ e_3e_1=\frac{3}{2}e_5$&$ e_3e_2=-3e_6  $\\&$ e_4e_1= \frac{1}{2}e_5 $
     &$ e_4e_2=e_6  $&$  e_5e_1=-4e_6  $& &\\
    \\
    ${\mathcal A}^6_{15}:  $&$ e_1e_1=e_2  $&$ e_1e_2=e_4  $&$ e_1e_3= e_6  $&$ e_1e_4=\frac{1}{2}e_5 $
     
     &$ e_1e_5=2e_6$\\&$ e_2e_1=e_3  $&$ e_2e_2=e_5  $&$ e_2e_3=-3e_6  $ 
     &$ e_2e_4=-2e_6  $&$ e_3e_1=\frac{3}{2}e_5$\\&$ e_3e_2=-3e_6  $&$ e_4e_1= \frac{1}{2}e_5+ e_6  $
     &$ e_4e_2=e_6  $&$  e_5e_1=-4e_6  $& &
\end{longtable}

  \subsection{Central extensions of ${\mathcal A}^{5}_{06}$}

Let us use the following notations: 

\[\nabla_1=[\Delta_{21}], 
\nabla_2=[\Delta_{15}]+[\Delta_{24}]+[\Delta_{33}]+[\Delta_{42}]+[\Delta_{51}].\] 
 
 The automorphism group of ${\mathcal A}^{5}_{06}$ consists of invertible matrices of the form

\[\phi= \left(\begin{array}{ccccc}
                               x & 0   & 0  & 0 & 0 \\
                               y & x^2 & 0  & 0 & 0 \\
                               z & 2xy & x^3 & 0 & 0 \\
                               v & 2xz + y^2 &3x^2y & x^4 & 0 \\
                               w & 2xv + 2yz & 3x^2z + 3xy^2 & 4x^3y & x^5
                               \end{array}\right)
                               .\]

Since
\[
\phi^T
                           \left(\begin{array}{ccccc}
                                0 & 0 & 0 & 0 & \alpha_2 \\
                                 \alpha_1 & 0 & 0 & \alpha_2 & 0 \\
                                 0 & 0 & \alpha_2 &  0 & 0\\
                            0 & \alpha_2 & 0 & 0 & 0\\
                            \alpha_2 & 0 & 0 & 0 & 0\\
                             \end{array}
                           \right)\phi
                           = \left(\begin{array}{ccccc}
                                \alpha^{****} & \alpha^{*} & \alpha^{**} & \alpha^{***} & \alpha^*_2 \\
                                 \alpha^{*}_1+\alpha^{*} & \alpha^{**} & \alpha^{***} & \alpha^*_2 & 0 \\
                                 \alpha^{**} & \alpha^{***} & \alpha^*_2 &  0 & 0\\
                            \alpha^{***} & \alpha^*_2 & 0 & 0 & 0\\
                            \alpha^*_2 & 0 & 0 & 0 & 0\\
                             \end{array}
                           \right),\]
we have that the action of ${\rm Aut} ({\mathcal A}^{5}_{06})$ on the subspace
$\langle  \sum\limits_{i=1}^2\alpha_i \nabla_i \rangle$
is given by
$\langle  \sum\limits_{i=1}^2\alpha^*_i \nabla_i \rangle,$ where

\[
\begin{array}{rclrcl}
\alpha^*_1&=&x^3 \alpha_1, &
\alpha^*_2&=& x^6 \alpha_2.\\
\end{array}\]

We suppose that $\alpha_2\neq 0$, otherwise obtained algebra gives an algebra with 2-dimensional annihilator. Therefore, consider the following cases:

\begin{enumerate}
\item If $\alpha_1=0,$ then $x=\frac{1}{\sqrt[6]{\alpha_2}},$  we have the representative $\langle \nabla_2\rangle;$
\item If $\alpha_1\neq0,$ then $x=\sqrt[3]{\frac{\alpha_1}{\alpha_2}},$  we have the representative $\langle \nabla_1+\nabla_2\rangle.$

\end{enumerate}

Hence, we have the following new algebras:

{\tiny
\begin{longtable}{lllllllllllllllllllllllll}

${\mathcal A}^{6}_{16}  $&$:$&$ 
e_1e_1=e_2 $ &$ e_1e_2=e_3  $&$ e_1e_3=e_4  $&$ e_1e_4=e_5  $&$ e_1e_5=e_6$&$ e_2e_1=e_3  $&$ e_2e_2=e_4 $ & $ e_2e_3=e_5$\\ && $e_2e_4=e_6$   &$ e_3e_1=e_4  $  &$ e_3e_2=e_5$ &$e_3e_3=e_6$&$ e_4e_1=e_5 $
&$e_4e_2=e_6$&$e_5e_1=e_6$\\

${\mathcal A}^{6}_{17}  $&$:$&$ 
e_1e_1=e_2 $ &$ e_1e_2=e_3  $&$ e_1e_3=e_4  $&$ e_1e_4=e_5  $&$ e_1e_5=e_6$&$ e_2e_1=e_3+e_6  $&$ e_2e_2=e_4 $ & $ e_2e_3=e_5$\\ && $e_2e_4=e_6$   &$ e_3e_1=e_4  $  &$ e_3e_2=e_5$ &$e_3e_3=e_6$&$ e_4e_1=e_5 $
&$e_4e_2=e_6$&$e_5e_1=e_6$\\
\end{longtable}}
 
  \subsection{Central extensions of ${\mathcal A}^{5}_{07}$}

Let us use the following notations: 
\[\nabla_1=[\Delta_{21}], \ \nabla_2=[\Delta_{15}]+2[\Delta_{22}]+[\Delta_{24}]+3[\Delta_{31}]+[\Delta_{33}]+[\Delta_{42}]+[\Delta_{51}].
\]
The automorphism group of ${\mathcal A}^{5}_{07}$ consists of invertible matrices of the form

\[\phi_i= \left(\begin{array}{ccccc}
                               (-1)^k & 0   & 0  & 0 & 0 \\
                               x & 1 & 0  & 0 & 0 \\
                               y & (-1)^k 2x & (-1)^k & 0 & 0 \\
                               z & x^2+(-1)^k2y &3x & 1 & 0 \\
                               t & 2xy+(-1)^k(x+2z) & (-1)^k3x^2+3y & (-1)^k4x & (-1)^k
                               \end{array}\right)
                               ,\]
where $k\in\{1,2\}.$ Since
\[
\phi_i^T
                           \left(\begin{array}{ccccc}
                                0 & 0 & 0 & 0 & \alpha_2 \\
                                 \alpha_1 & 2\alpha_2 & 0 & \alpha_2 & 0 \\
                                 3\alpha_2 & 0 & \alpha_2 &  0 & 0\\
                            0 & \alpha_2 & 0 & 0 & 0\\
                            \alpha_2 & 0 & 0 & 0 & 0\\
                             \end{array}
                           \right)\phi_i
                           = \left(\begin{array}{ccccc}
                                \alpha^{****} & \alpha^{*} & \alpha^{**} & \alpha^{***} & \alpha^*_2 \\
                                 \alpha^*_1+\alpha^{*} & 2\alpha^*_2+\alpha^{**} & \alpha^{***} & \alpha^*_2 & 0 \\
                                 3\alpha^*_2 & \alpha^{***} & \alpha^*_2 &  0 & 0\\
                            \alpha^{***} & \alpha^*_2 & 0 & 0 & 0\\
                            \alpha^*_2 & 0 & 0 & 0 & 0\\
                             \end{array}
                           \right),\]
we have that the action of ${\rm Aut} ({\mathcal A}^{5}_{07})$ on the subspace
$\langle  \sum\limits_{i=1}^2\alpha_i \nabla_i \rangle$
is given by
$\langle  \sum\limits_{i=1}^2\alpha^*_i \nabla_i \rangle,$ where

\[
\begin{array}{rclrcl}
\alpha^*_1&=&(-1)^i \alpha_1-6x\alpha_2, &
\alpha^*_2&=& \alpha_2.\\
\end{array}\]

We have only one non-trivial orbit with the representative $\langle \nabla_2 \rangle,$ and get 

 \begin{longtable}{lllllllllllllllllllllllll}
  ${\mathcal A}^{6}_{18}  $&$:$&$ e_1e_1=e_2  $&$ e_1e_2=e_3  $&$ e_1e_3=e_4  $&$ e_1e_4=e_5  $ &$ e_1e_5=e_6$\\ &&$ e_2e_1=e_3+e_5  $&$ e_2e_2=e_4+2e_6   $&$ e_2e_3=e_5  $&$ e_2e_4=e_6  $&$  e_3e_1=e_4+3e_6$\\
  &&$  e_3e_2=e_5  $&$ e_3e_3=e_6$&$ e_4e_1=e_5  $&$ e_4e_2=e_6  $&$ e_5e_1=e_6 $\\
\end{longtable}

\subsection{Classification theorem}
Summarizing results of the present section we have the following theorem.

\begin{theoremB}
Let $\mathcal A$ be a $6$-dimensional complex one-generated nilpotent assosymmetric  algebra,
then $\mathcal A$ is isomorphic to an algebra from the following list.	

{\tiny
\begin{longtable}{|ll lll|}

\hline

    ${\mathcal A}^6_{01}$&$:$&$ 
    e_1e_1=e_2 $&$ e_1e_2=e_4$&$ e_1e_3=e_5$\\
    && $e_1e_4=e_6$ &$e_2e_1=e_3 $ &$e_3e_1=-e_6 $\\
    &&$  e_4e_1=e_5-e_6$&&\\
\hline
     
    ${\mathcal A}^6_{02}(\alpha) $&$:$&$ 
    e_1e_1=e_2 $&$ e_1e_2=e_4 $&$ e_1e_3=e_5 $\\\ 
    &&$ e_1e_4=\alpha e_6 $ &$  e_2e_1=e_3 $ &$ e_2e_2=e_6 $\\
    &&$  e_3e_1=(2-\alpha) e_6$& $ e_4e_1=e_5-(\alpha-1)e_6$&\\
     
\hline
    ${\mathcal A}^6_{03}(\alpha,\beta) $&$:$&
    $ e_1e_1=e_2 $&$ e_1e_2=e_4 $&$ e_1e_3=\alpha e_5+\beta e_6 $\\
    &&$ e_1e_4=e_5 $ &$ e_2e_1=e_3$&$ e_2e_2=e_6 $ \\
    &&$ e_3e_1=- e_5+2e_6 $&  $e_4e_1=(\alpha-1) e_5+(\beta+1)e_6$ & \\
    
\hline
     ${\mathcal A}^6_{04}$ &$ : $&
     $ e_1e_1=e_2 $&$ e_1e_2=e_3 $&$ e_1e_3=e_4 $\\
     &&$ e_1e_4=e_5 $&$ e_2e_1=e_3+e_6$ &$ e_2e_2=e_4 $\\
     &&$ e_2e_3=e_5 $&$ e_3e_1=e_4 $&$    e_3e_2=e_5 $\\
     &&$ e_4e_1=e_5$ &&\\

\hline
  $ {\mathcal A}^{6}_{05}(\alpha)$ &$: $&$  
  e_1e_1=e_2  $&$ e_1e_2=e_3   $&$  e_1e_3=-e_5+\alpha e_6  $\\
  &&$ e_1e_4= e_6  $ &$ e_1e_5=e_6$&$ e_2e_1=e_3 +e_4$ \\
   &&$ e_2e_2=-e_5   $&$ e_2e_3=-e_6 $&$ e_3e_1=-e_5+ e_6  $\\
       &&$ e_3e_2=-e_6$& $ e_4e_1=-(\alpha+1)  e_6  $&$ e_5e_1=e_6 $\\
   
\hline
${\mathcal A}^{6}_{06}$  &$: $&
$e_1e_1=e_2  $&$ e_1e_2=e_3   $&$  e_1e_3=-e_5+ e_6  $\\
&&$ e_1e_5=e_6  $ &$ e_2e_1=e_3 +e_4$  &$ e_2e_2=-e_5  $\\
&&$ e_2e_3=-e_6 $&$ e_3e_1=-e_5  $&$ e_3e_2=-e_6  $\\
&&$ e_4e_1=- e_6  $&$ e_5e_1=e_6 $&\\
   
\hline
$ {\mathcal A}^{6}_{07}$ &$: $&
$  e_1e_1=e_2  $&$ e_1e_2=e_3   $&$  e_1e_3=-e_5  $\\
&&$ e_1e_5=e_6  $  &$ e_2e_1=e_3 +e_4$ & $ e_2e_2=-e_5$\\
&&$ e_2e_3=-e_6 $&$ e_3e_1=-e_5  $&$ e_3e_2=-e_6  $\\
&&$ e_5e_1=e_6 $ && \\   
\hline

\multicolumn{5}{|l|}{
$i=08$\mbox{ for }$\alpha=\frac{1}{2}(\beta +\sqrt{-2+6\beta -3 \beta ^2})$
\mbox{ and  }
$i=10$\mbox{ for }$\alpha=\frac{1}{2}(\beta-\sqrt{-2+6\beta -3 \beta ^2})$ \mbox{ with } $\beta\neq\frac{1}{2}$}\\
 ${\mathcal A}^6_{i}(\beta) $&$:$& 
    $e_1 e_1=e_2 $&$  e_1 e_2=e_4$ & $ e_1 e_3=\alpha e_5$ \\
    &&$  e_1 e_4=\beta e_5 $ &$ e_1 e_5=(2 \beta-1) e_6 $& $ e_2 e_1=e_3$ \\
    &&$  e_2 e_2=e_5 $&$  e_2 e_3= (2 \alpha \beta-2 \beta+1)e_6  $ &$  e_2 e_4= (\alpha+2 \beta^2-3 \beta+1)e_6 $\\ 
    &&$ e_3 e_1=(2-\beta) e_5 $&$  e_3 e_2=(3 \alpha-2 \alpha \beta+2 \beta^2-3 \beta+1)e_6 $&$   e_4 e_1= (\alpha-\beta+1)e_5$
    \\&&$ e_4 e_2=( 2 \alpha-2 \alpha \beta+2 \beta-1)e_6$  &$  e_5 e_1= (2 \alpha-2 \beta+1)e_6$ &\\

  \hline
  \multicolumn{5}{|l|}{
$i=09$\mbox{ for }$\alpha=\frac{1}{2}(\beta +\sqrt{-2+6\beta -3 \beta ^2})$ \mbox{ with } $\beta\neq\frac{3}{2}$
\mbox{ and  }
$i=11$\mbox{ for }$\alpha=\frac{1}{2}(\beta-\sqrt{-2+6\beta -3 \beta ^2})$ \mbox{ with } $\beta\neq\frac{1}{2}$ }\\

    ${\mathcal A}^6_{i+1}(\beta) $&$:$& 
    $e_1 e_1=e_2 $&$  e_1 e_2=e_4$ & $ e_1 e_3=\alpha e_5$ \\
    &&$  e_1 e_4=\beta e_5+e_6 $ &$ e_1 e_5=(2 \beta-1) e_6 $& $ e_2 e_1=e_3$ \\
    &&$  e_2 e_2=e_5 $&$  e_2 e_3= (2 \alpha \beta-2 \beta+1)e_6  $ &$  e_2 e_4= (\alpha+2 \beta^2-3 \beta+1)e_6 $\\ 
    &&$ e_3 e_1=(2-\beta) e_5-e_6 $&$  e_3 e_2=( 3 \alpha-2 \alpha \beta+2 \beta^2-3 \beta+1)e_6 $&$   e_4 e_1= (\alpha-\beta+1)e_5-e_6$
    \\&&$ e_4 e_2=(2 \alpha-2 \alpha \beta+2 \beta-1)e_6$  &$  e_5 e_1= (2 \alpha-2 \beta+1)e_6$ &\\    
    \hline
${\mathcal A}^6_{12}(\alpha)$& $:$ & 
$e_1e_1=e_2$ & $e_1e_2=e_4$  & $e_1e_3= e_5+ e_6$ \\
&& $e_1e_4= e_5+2\alpha e_6 $ &  $e_1e_5=e_6$ &$e_2e_1=e_3$ \\
&& $e_2e_2=e_5+\alpha e_6 $& $e_2e_3=e_6$ & $e_2e_4=e_6$ \\
&& $e_3e_1=e_5$ &  $e_3e_2=e_6$ & $e_4e_1=e_5+(1-\alpha) e_6$\\
&& $e_4e_2=e_6$ & $e_5e_1=e_6$ &\\
    \hline   

${\mathcal A}^6_{13}$&$:$ &$ e_1e_1=e_2  $&$ e_1e_2=e_4  $&$ e_1e_3= e_5+e_6 $\\&&$ e_1e_4=\frac{3}{2}e_5$&$ e_1e_5=2 e_6$
     &$ e_2e_1=e_3$ \\&&$ e_2e_2=e_5  $&$e_2e_3=e_6  $&$ e_2e_4=2e_6  $   \\&&$  e_3e_1=\frac{1}{2}e_5  $   &$ e_3e_2=e_6 $ &$ e_4e_1=e_6  $\\& &$ e_4e_2=e_6  $   &     &  \\     
       \hline
     ${\mathcal A}^6_{14}$&$:  $&$ e_1e_1=e_2  $&$ e_1e_2=e_4  $&$ e_1e_4=\frac{1}{2}e_5 $\\
     
     &&$ e_1e_5=2e_6$&$ e_2e_1=e_3  $&$ e_2e_2=e_5  $\\ &&$ e_2e_3=-3e_6  $ 
     &$ e_2e_4=-2e_6  $&$ e_3e_1=\frac{3}{2}e_5$\\
     &&$ e_3e_2=-3e_6  $&$ e_4e_1= \frac{1}{2}e_5 $
     &$ e_4e_2=e_6  $\\&&$  e_5e_1=-4e_6  $& &\\
    \hline
    ${\mathcal A}^6_{15}$&$:  $&$ e_1e_1=e_2  $&$ e_1e_2=e_4  $&$ e_1e_3= e_6  $\\
    &&$ e_1e_4=\frac{1}{2}e_5 $
     
     &$ e_1e_5=2e_6$&$ e_2e_1=e_3  $\\
     &&$ e_2e_2=e_5  $&$ e_2e_3=-3e_6  $ 
     &$ e_2e_4=-2e_6  $\\
     &&$ e_3e_1=\frac{3}{2}e_5$&$ e_3e_2=-3e_6  $&$ e_4e_1= \frac{1}{2}e_5+ e_6  $\\
     &&$ e_4e_2=e_6  $&$  e_5e_1=-4e_6  $&\\
       
       \hline
${\mathcal A}^{6}_{16}  $&$:$&$ 
e_1e_1=e_2 $ &$ e_1e_2=e_3  $&$ e_1e_3=e_4  $\\
&&$ e_1e_4=e_5  $&$ e_1e_5=e_6$&$ e_2e_1=e_3  $\\
&&$ e_2e_2=e_4 $ & $ e_2e_3=e_5$ & $e_2e_4=e_6$\\   
&&$ e_3e_1=e_4  $  &$ e_3e_2=e_5$ &$e_3e_3=e_6$\\
&&$ e_4e_1=e_5 $&$e_4e_2=e_6$&$e_5e_1=e_6$\\
\hline

${\mathcal A}^{6}_{17}  $&$:$&$ 
e_1e_1=e_2 $ &$ e_1e_2=e_3  $&$ e_1e_3=e_4  $\\
&&$ e_1e_4=e_5  $&$ e_1e_5=e_6$&$ e_2e_1=e_3+e_6  $\\
&&$ e_2e_2=e_4 $ & $ e_2e_3=e_5$ & $e_2e_4=e_6$ \\  
&&$ e_3e_1=e_4  $  &$ e_3e_2=e_5$ &$e_3e_3=e_6$\\
&&$ e_4e_1=e_5 $&$e_4e_2=e_6$&$e_5e_1=e_6$\\
\hline

  ${\mathcal A}^{6}_{18}  $&$:$&
  $ e_1e_1=e_2  $&$ e_1e_2=e_3  $&$ e_1e_3=e_4  $\\
  &&$ e_1e_4=e_5  $ &$ e_1e_5=e_6$ &$ e_2e_1=e_3+e_5  $\\
  &&$ e_2e_2=e_4+2e_6   $&$ e_2e_3=e_5  $&$ e_2e_4=e_6  $\\
  &&$ e_3e_1=e_4+3e_6$ &$  e_3e_2=e_5  $&$ e_3e_3=e_6  $\\ 
  &&$ e_4e_1=e_5  $&$ e_4e_2=e_6  $&$ e_5e_1=e_6 $\\
\hline
\end{longtable}
}
\begin{center}
    
Note: ${\mathcal A}^6_{08}(\frac{3}{2})\cong {\mathcal A}^6_{09}(\frac{3}{2}).$
\end{center}

\end{theoremB}

\end{document}